\documentclass[11pt]{article}
\usepackage[table]{xcolor} 
\usepackage{amsmath,amsfonts,amssymb}
\usepackage{graphicx}
\usepackage{setspace}
\usepackage[flushleft]{threeparttable}
\usepackage{subcaption}
\usepackage{epsfig}
\usepackage{color}
\usepackage{multirow}
\usepackage{amsmath,amsfonts,amssymb}
\usepackage{graphicx}
\usepackage{setspace}
\usepackage{subcaption}
\usepackage{epsfig}
\usepackage{multirow}
\usepackage{booktabs}
\usepackage{lscape}
\usepackage{enumerate}
\usepackage{enumitem}
\usepackage{mathtools}
\usepackage{relsize}
\usepackage{slantsc}
\usepackage[square,numbers]{natbib}
\usepackage[top=1in, bottom=1in, left=1in, right=1in]{geometry}
\usepackage[utf8]{inputenc}
\usepackage{tikz}
\usepackage[ruled,linesnumbered]{algorithm2e}
\usetikzlibrary{positioning}

\usepackage{authblk}

\usepackage{lmodern}
 \tikzset{
    invisible/.style={opacity=0,text opacity=0},
    visible on/.style={alt={#1{}{invisible}}},
    alt/.code args={<#1>#2#3}{%
      \alt<#1>{\pgfkeysalso{#2}}{\pgfkeysalso{#3}} 
    },
  }
\usetikzlibrary{patterns}
\usepackage{pgfplots}
\pgfplotsset{compat=newest}
\usetikzlibrary{arrows,decorations.pathreplacing,backgrounds,automata}

\onehalfspace

\usepackage{ragged2e}

\newtheorem{theorem}{Theorem}

\newtheorem{proposition}[theorem]{Proposition}
\newtheorem{corollary}[theorem]{Corollary}
\newtheorem{example}{Example}

\newenvironment{proofin}[1][Proof:]{\begin{trivlist}\item[\hskip \labelsep {\bfseries #1}]}{\end{trivlist}}
\newenvironment{proof}[2] {\paragraph{Proof of {#1} {#2}:}}{\hfill$\square$}

\newcommand{\bA}{ \mathbf{A} }
\newcommand{\ba}{ \mathbf{a} }
\newcommand{\bbb}{ \mathbf{b} } 

\newcommand{\bc}{ \mathbf{c} }

\newcommand{\beee}{ \mathbf{e} }

\newcommand{\bu}{ \mathbf{u} }

\newcommand{\bv}{ \mathbf{v} }

\newcommand{\bx}{ \mathbf{x} }

\newcommand{\by}{ \mathbf{y} }

\newcommand{\bz}{ \mathbf{z} }
\newcommand{\bbR}{ \mathbb{R} }
\newcommand{\cI}{ \mathcal{I} }

\newcommand{\cJ}{ \mathcal{J} }
\newcommand{\cX}{ \mathcal{X} }
\newcommand{\cK}{ \mathcal{K} }
\newcommand{\cL}{ \mathcal{L} }

\newcommand{\cQ}{ \mathcal{Q} }

\newcommand{\cG}{ \mathcal{G} }

\newcommand{\bxhat}{\mathbf{\hat{x}}}

\newcommand{\bzero}{ \mathbf{0} }

\newcommand{\blambda}{ \boldsymbol{\lambda} }

\newcommand{\cen}{\text{cen}}

\def\argmin{\mathop{\rm argmin}}
\def\argmax{\mathop{\rm argmax}}

\def\zahed#1{{\color{black}#1}} 

\def\rev#1{{\color{black}#1}}

\renewenvironment{abstract}
{\begin{quote}
\noindent \rule{\linewidth}{.5pt}\par{\bfseries \abstractname.}}
{\medskip\noindent \rule{\linewidth}{.5pt}
\end{quote}
}
\begin{document}
\title{Optimality-Based Clustering: An Inverse Optimization Approach}

\author[a]{\small Zahed Shahmoradi}
\author[a]{\small Taewoo Lee\thanks{Corresponding Author. tlee6@uh.edu}}

\affil[a]{\footnotesize Department of Industrial Engineering, University of Houston, Houston, TX 77204, USA}

\date{}
\maketitle
\vspace{-0.8in}
\begin{abstract}
We propose a new clustering approach, called optimality-based clustering, that clusters data points based on their latent decision-making preferences. We assume that each data point is a decision generated by a decision-maker who (approximately) solves an optimization problem and cluster the data points by identifying a common objective function of the optimization problems for each cluster such that the worst-case optimality error is minimized. We propose three different clustering models \rev{and 
test them in the diet recommendation application}.
\begin{flushleft}
\end{flushleft}
\vspace{-0.7cm}
\end{abstract}

\vspace{-0.2in}
\section{Introduction}\label{sec:Clust:Introduction}
\vspace{-0.1in}
Clustering is a technique that groups objects (e.g., data points) into clusters such that the objects in the same cluster are more similar to each another than to those in other clusters based on some similarity measure \citep{hastie2009elements}. Most clustering approaches fall into one of the following categories: centroid-based clustering, distribution-based clustering, and density-based clustering. In centroid-based clustering, each object is assigned to a cluster based on its similarity to a representative object called a centroid (e.g., K-means clustering) \citep{jain1999data,likas2003global,xu2008clustering}. In density-based clustering, a density measure, e.g., the number of objects within a certain distance, is used to detect areas with high density, in which the objects are grouped into the same cluster  \citep{ester1996density,kriegel2011density}. Distribution-based clustering groups the objects based on whether or not they belong to the same distribution \citep{xu1998distribution}.

Often, data points correspond to decisions generated by decision-makers (DMs) who are assumed to solve some kind of decision-making problems (DMPs). Although traditional clustering approaches for such decision data may indicate which decisions are similar to each other, this similarity does not necessarily imply that the DMs whose decisions are in the same cluster have similar preferences. For example, suppose the DMPs can be formulated as optimization problems where the DM's preferences are encoded in the objective function parameters. Even when two DMs' decisions are geometrically close to each other, they might have been generated by two DMPs with completely different objective function parameters under different feasible regions, which traditional clustering cannot capture. The focus of this paper is to cluster decision data based on the similarity in the DM's decision-making preferences, captured by parameters in their underlying DMPs.

\rev{
Clustering based on decision-making preferences can help create targeted, group-based decision support tools. For example, by clustering patients based on their health-related preferences (e.g., health benefit vs. cost saving) using their past disease screening decisions, one can create a group-based yet easily implementable screening guideline that is consistent with the patients' preferences (e.g., increased use of telemedicine for a specific group of patients). Similarly, when developing a diet recommendation system, clustering individuals based on their food preferences and inferring a common objective function for each cluster can help create a group-specific diet recommendation framework. A post-hoc analysis can be done to further identify association of the preference clusters with other factors such as health conditions and socio-demographic factors.
}

Since this clustering problem requires inferring objective function parameters of the DMPs from decision data, it inherently involves inverse optimization. Given an observed decision from a DM who solves an optimization problem, inverse optimization infers parameters of the problem that make the decision as optimal as possible (e.g., \citep{ahuja2001inverse, aswani2018inverse,bertsimas2015,chan2014,esfahani2018data,keshavarz2011imputing}).  Solving the DM's optimization problem with these inferred parameters then leads to a decision that is close to the observed one. Previous inverse optimization models assume that decision data is collected from either a single DM or a group of DMs whose preferences are known to be similar, for which the same, single set of parameters is inferred \citep{aswani2018inverse,babier2021ensemble,esfahani2018data,shahmoradi2019quantile}.

In this paper, we develop a new clustering approach that clusters decision data (hence DMs) based on their latent decision-making preferences. In particular, inspired by inverse optimization, we propose the clustering problem that simultaneously groups observed decisions into clusters and finds an objective function for each cluster such that the decisions in the same cluster are rendered as optimal as possible for the assumed DMPs. We 
use optimality errors associated with the decisions with respect to the inferred objective function as a measure of similarity; hence we call this problem ``optimality-based clustering." We further enhance the problem by incorporating the notion of cluster stability, measured for each cluster by the worst-case distance between the decision data in the cluster and optimal decisions achieved by the DMPs using the inferred objective function for the cluster. The stability-driven, optimality-based clustering problem is computationally challenging. We derive mixed-integer programs (MIPs) that provide upper and lower bound solutions for the true clustering problem as well as heuristics that approximately solve this problem. \rev{Finally, we demonstrate the proposed clustering approach in the diet recommendation application to cluster individuals based on their food preferences. Unless otherwise stated, proofs are in the appendix.}

\vspace{-0.1in}
\section{Preliminaries}\label{sec:Clust:Prelim}
\vspace{-0.1in}
In this section, we present an initial formulation for the optimality-based clustering problem and a simple example to demonstrate the idea. We then define the notion of cluster stability in the context of optimality-based clustering, which we later use to propose an enhanced clustering formulation.
\vspace{-0.1in}
\subsection{A General Clustering Problem}
\vspace{-0.05in}
We focus on a centroid-based clustering problem where the similarity of a data point to a cluster is assessed by the distance between the data point and a centroid of the cluster. Given a dataset $\hat\cX=\{\bxhat^1,\ldots,\bxhat^K\}$ with the index set $\cK=\{1,\ldots,K\}$, let $\{\cG^\ell\}_{\ell\in \cL}$ be a collection of $L$ clusters where $\cG^\ell\subseteq\cK$ and $\cL=\{1,\ldots,L\}$. For each cluster $\cG^\ell$, the {\it{dissimilarity}} among the members of the cluster is measured by 
$\displaystyle\sum_{k\in \cG^\ell}d(\bxhat^k,\bx^\ell),$ 
where $\bx^\ell$ denotes the centroid of the cluster and $d(\bxhat^k,\bx^\ell)$ represents the distance between observation $\bxhat^k$ and its cluster centroid $\bx^\ell$, e.g., $d(\bxhat^k,\bx^\ell)=\lVert \bxhat^k-\bx^\ell \rVert_r$ for some $r\ge 1$. Based on the above definition, a centroid-based clustering problem seeks clusters $\{\cG^\ell\}_{\ell\in \cL}$ such that the sum of dissimilarities over all clusters is minimized, i.e.,  
\begin{equation}\label{eq:Cluster}\nonumber
\begin{aligned}
 \underset{\{\cG^\ell\}_{\ell\in \cL}, \{\bx^\ell\}_{\ell\in \cL}}{\text{minimize}} &\quad  \sum_{\ell\in \cL}  \sum_{k\in \cG^\ell}d(\bxhat^k,\bx^\ell).\\
\end{aligned}
\end{equation}

\subsection{Optimality-Based Clustering: The Initial Model}\label{sec:Clust:Example}
We assume that each data point $\hat\bx^k\in \hat\cX$ is an observed decision from DM $k$ (denoted by $\text{DM}_k$) who approximately solves the following optimization problem as a decision-making problem ($\text{DMP}_k$):
\begin{equation}\label{eq:FO}\nonumber
\textrm{\textbf{DMP}}_k(\bc):
\quad \underset{\bx}{\text{minimize}}\;\{{\bc}'\bx\,|\, \bA^k\bx \geq \bbb^k\},
\end{equation}
where $\bc \in \bbR^n, \bx \in \bbR^n$, $\bA^k \in \bbR^{m_k\times n}$, and $\bbb^k \in \bbR^{m_k}$, for each $k\in \cK$. For each DM $k\in \cK$, let $\cI^k=\{1,\ldots,m_k\}$ and $\cJ=\{1,\ldots,n\}$ index the constraints and variables of DMP${_k}$, and $\ba^{ki} \in\mathbb{R}^n$ be a (column) vector corresponding to the $i$-th row of $\bA^k$. We let $\cX^k$ be the set of feasible solutions for $\text{DMP}_k$, assumed bounded,  full-dimensional, and free of redundant constraints, and $\cX^{ki}=\{\bx\in \cX^k \,|\, {\ba^{ki}}'\bx=b_i^k\}$, $i\in\cI^k$. Let $\cX^{k*}(\bc) = \argmin \textbf{DMP}_k(\bc)$. Without loss of generality, we assume that each $\ba^{ki}$ 
is normalized {\it{a priori}} 
such that $\|\ba^{ki}\|_1=1$. 

\begin{figure}[h]\centering
\centering
\vfill
\begin{tikzpicture}[scale=1.0]
\draw[->] (0,0) -- coordinate (x axis mid) (3,0) node[pos=1,below=0.1cm] {$x_1$};
    \foreach \x in {0,0.5,...,2.5}
		\draw  (\x,1pt)--(\x,-1pt);
    \foreach \x in {0,...,2.5}
		\draw  (\x,1pt)--(\x,-1pt)
		node[anchor=north] {\x};
    \draw[->] (0,0) -- coordinate (y axis mid) (0,3) node[pos=1,left] {$x_2$};
    \foreach \y in {0,0.5,...,2.5}
		\draw  (1pt,\y)--(-1pt,\y);
    \foreach \y in {0,...,2.5}
		\draw  (1pt,\y)--(-1pt,\y)
		node[anchor=east] {\y};
		

	\draw (1.5,0) -- (1.5,1);
	\draw[gray,dotted] (1.5,0)--(1.5,-0.5);
	\draw[gray,dotted] (1.5,1)--(2,1);

	\draw (1.5,1) -- (0,1);
	\draw[gray,dotted] (-0.5,1) -- (0,1);
	\draw[gray,dotted] (1.5,1) -- (1.5,1.5);

    \draw[fill=black] (1.2,1) circle (0.04);
    \node [black,anchor= south] at (1.2,1) {\scriptsize$\hat\bx^1$};
    
    \draw[fill=black] (1.5,0.6) circle (0.04);
    \node [black,anchor= east] at (1.5,0.6) {\scriptsize$\hat\bx^2$};
    
     \node [black,anchor= west] at (1.5,0.3) {
     };
    \node [black,anchor= south] at (1,1) 
    {
    };

    \draw (2.5,0) -- (2.5,2.5);
	\draw[gray,dotted] (2.5,0)--(2.5,-0.5);
	\draw[gray,dotted] (2.5,2.5)--(2.5,3);

	\draw (2.5,2.5) -- (0,2.5);
	\draw[gray,dotted] (-0.5,2.5) -- (0,2.5);
	\draw[gray,dotted] (2.5,2.5) -- (3,2.5);

    \draw[fill=black] (2.5,0.3) circle (0.04);
    \node [black,anchor= east] at (2.5,0.3) {\scriptsize$\hat\bx^3$};
    
    \node [black,anchor= west] at (2.5,1.25) {
    };
    \node [black,anchor= south] at (1.25,2.5) {
    };
    
\end{tikzpicture}
\vfill
\caption{
Observations from DMs $k=1,2,$ and $3$ and their respective feasible regions. 
}
\label{fig:Clustering_Motiv}
\end{figure}
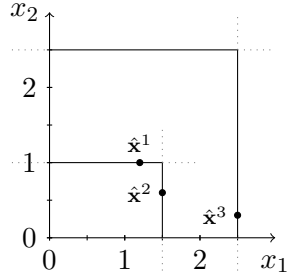

Given a set of observed decisions $\hat\cX$, the goal of optimality-based clustering is to group the observations into $L < K$ clusters $\{\cG^\ell\}_{\ell\in \cL}$ and find a cost vector $\bc^\ell$ for each cluster $\ell$ such that each observation $\hat\bx^k$ in cluster $\ell$ (i.e., for $k\in \cG^\ell$) is as close as possible to an optimal solution to $\text{DMP}_k(\bc^{\ell})$
. This problem can be formulated as follows:   
\begin{subequations}\label{eq:Abstract_OptClust}
\begin{align}
\quad \underset{\{\bx^k\}_{k\in \cK},\{(\bc^\ell,\cG^\ell)\}_{\ell\in \cL}}{\text{minimize}}  &\quad  d(\hat\cX,\{\bx^k\}_{k\in \cK})  \label{eq:Abstract_OptClust_1}\\
 \text{subject to} &\quad \bx^k \in 
\cX^{k*}(\bc^{\ell}),  \quad \forall k\in \cG^\ell,\ell \in \cL, \label{eq:Abstract_OptClust_2}\\
& \quad \lVert\bc^\ell\rVert_1=1, \quad \forall \ell\in \cL. \label{eq:Abstract_OptClust_3} 
\end{align}
\end{subequations}
The objective of the above problem is to minimize the distance between the observations $\hat\cX$ and solutions $\bx^k$'s that are optimal for their respective DMPs with respect to $\bc^{\ell}$ for $k\in\cG^\ell$; e.g., $\displaystyle d(\hat\cX,\{\bx^k\}_{k\in \cK})=\sum_{\ell\in \cL}  \sum_{k\in \cG^\ell}\|\bxhat^k-\bx^k\|$. Constraint \eqref{eq:Abstract_OptClust_3} prevents the trivial solution $\bc^\ell=\bzero$ from being feasible
. Note that the above problem is analogous to centroid-based clustering problems in that $\bc^\ell$ can be seen as the centroid of cluster $\ell$, representing the shared decision preference of the observations assigned to cluster $\ell$. 
We use the following simple example to demonstrate the idea.
%
\begin{example}\label{example:clustering}
Suppose three DMs solve the following problem with their own objective functions: 
\begin{equation}\nonumber
\underset{x_1,x_2}{\textup{maximize}}\;\{c_1 x_1 + c_2 x_2 \,|\, x_1\le b_1,\,x_2 \le b_2,\,x_1,x_2 \ge 0\}.
\end{equation}
%
Let $(b_1,b_2)=(1.5,1)$ for DMs $k=1,2$ and $(b_1,b_2)=(2.5,2.5)$ for DM $k=3$ (see Figure \ref{fig:Clustering_Motiv} for the feasible regions). We assume the following decisions are observed from the DMs: $\hat\bx^1=\begin{bsmallmatrix}1.2\\1\end{bsmallmatrix}$, $\hat\bx^2=\begin{bsmallmatrix}1.5\\0.6\end{bsmallmatrix}$, and $\hat\bx^3=\begin{bsmallmatrix}2.5\\0.3\end{bsmallmatrix}$ (see Figure \ref{fig:Clustering_Motiv}). If the desired number of clusters is two (i.e., $L=2$), traditional K-means clustering based on the Euclidean distance finds 
$\{\hat\bx^1,\hat\bx^2\}$ and 
$\{\hat\bx^3\}$ to be optimal clusters. 
However, if the goal is to group the decisions based on the preferences encoded in the corresponding DMPs, clustering should be done differently. In particular, given their respective feasible regions, $\hat\bx^2$ and $\hat\bx^3$ share the same preference as they are optimal for their respective DMPs based on the same cost vector 
$\bc=\begin{bsmallmatrix}1\\0\end{bsmallmatrix}$; on the other hand, $\hat\bx^1$ is optimal to the DMP with respect to a completely different cost vector 
$\bc=\begin{bsmallmatrix}0\\1\end{bsmallmatrix}$. As a result, 
an optimal clustering is $\{\hat\bx^2,\hat\bx^3\}$ and $\{\hat\bx^1\}$.
\end{example}
\vspace{-0.14in}
\vspace{-0.2in}
\subsection{Cluster Instability}
\vspace{-0.03in}
In this subsection, we show that the initial model~\eqref{eq:Abstract_OptClust} 
is often subject to an instability issue due to the structure of the DMP formulation and propose a measure of instability in the context of optimality-based clustering. Given an optimal cost vector $\bc^{\ell*}$ for some cluster $\ell$ achieved by model~\eqref{eq:Abstract_OptClust}, we note that \textbf{DMP}$_k(\bc^{\ell*})$ often leads to $\bx^k\in
\cX^{k*}(\bc^{\ell})$ that is far from the observations assigned to the cluster. 
For illustration, consider the same example in Figure \ref{fig:Clustering_Motiv}, where model \eqref{eq:Abstract_OptClust} finds $\cG_2=\{\hat\bx^{1}\}$ (i.e., $\hat\bx^1$ assigned to cluster $\ell=2$) and $\bc^{2*}=\begin{bsmallmatrix}0\\1\end{bsmallmatrix}$. While the desirable forward optimal solution with respect to this cost vector is supposed to be close to $\hat\bx^1$, solving $\text{DMP}_1(\bc^{2*})$ can lead to an optimal solution $\bx^*=\begin{bsmallmatrix}0\\1\end{bsmallmatrix}$, which is far from $\hat\bx^1$. Note that this cluster instability issue is different from the cluster assignment instability issues considered in the traditional clustering literature \citep{rakhlin2007stability,von2010clustering}; it is rather associated with the argmin set of the DMP for a certain cost vector. This type of instability is also discussed in Shahmoradi and Lee \cite{shahmoradi2019quantile} in the context of inverse linear programming.

We now formally define a notion of \textit{cluster stability}, which we then use to propose an enhanced clustering problem formulation that 
improves on the initial model \eqref{eq:Abstract_OptClust} in the next section. 
%
%
Given that the instability issue is caused by $\bx\in\cX^{k*}(\bc)$ being too far from $\hat\bx^{k}$, we assess the instability of a cost vector $\bc$ associated with each $\hat\bx^{k}$ via the worst-case distance between $\hat\bx^{k}$ and $\cX^{k*}(\bc)$:
\begin{equation}\label{eq:per_k_stability}\nonumber
{\max}\, \{d(\hat\bx^k,\bx^{k}) \,|\, \bx^{k} \in \cX^{k*}(\bc)\}.
\end{equation}
%
Then, the instability of cluster $\cG^\ell$ with its cost vector $\bc^{\ell}$ is measured by the following measure:
\begin{equation}\label{eq:cluster_stability}
\, \underset{k\in \cG^\ell}{\max} \, \max\, \{d(\hat\bx^k,\bx^{k}) \,|\, \bx^{k} \in \cX^{k*}(\bc^\ell)\}.
\end{equation}
In other words, cluster $\ell$ is said to be more stable if its cost vector $\bc^\ell$ leads to a smaller worst-case distance between $\hat\bx^k$ and the set of optimal solutions for DMP$_{k}(\bc^{\ell})$ over all data points in the cluster. For brevity, from here on out we combine the two max terms in \eqref{eq:cluster_stability} and simply write it as $\underset{k\in \cG^\ell}{\max} \, \{d(\hat\bx^k,\bx^{k}) \,|\, \bx^{k} \in \cX^{k*}(\bc^\ell)\}.$

%



\vspace{-0.05in}
\section{Models}\label{sec:Clust:Models}
\vspace{-0.1in}
In this section, we first propose an enhanced optimality-based clustering problem that addresses the cluster instability issue by incorporating the stability measure \eqref{eq:cluster_stability}. We also propose two heuristics that approximately solve the problem by separating it into two stages: the clustering stage and the cost vector inference stage. 
We then analytically compare the performances of these approaches.
%
\vspace{-0.1in}
\subsection{The Stability-Driven Clustering Model}\label{sec:Clust:Model:Stability}

To address the cluster instability issue in model~\eqref{eq:Abstract_OptClust}, we replace its objective function with the stability-incorporated dissimilarity measure in \eqref{eq:cluster_stability}. This leads to the following, which we call the stability-driven clustering (SC) model:
\begin{subequations}\label{eq:Abstract_StabClust}
\begin{align}
\textbf{SC}(\cK,L): \quad  \underset{\{(\bc^\ell,\cG^\ell)\}_{\ell\in \cL}}{\text{minimize}} & \quad \underset{\ell\in \cL}\max\; \; \underset{k\in \cG^\ell}{\max}\;\; \{d(\hat\bx^k,\bx^{k})\}  \label{eq:Abstract_StabClust_1}
\\
\text{subject to}  & \quad \bx^k \in 
\cX^{k*}(\bc^{\ell}),  \quad \forall \ell \in \cL, k\in \cG^\ell, \label{eq:Abstract_StabClust_2} \\
& \quad \lVert\bc^\ell\rVert_1=1, \quad \forall \ell\in \cL,\label{eq:Abstract_StabClust_3} 
\end{align}
\end{subequations}
where now the objective is to maximize stability for all clusters by minimizing the worst-case distance between $\hat\bx^k$ and the argmin set $\cX^{k*}(\bc^\ell)$ over all observations and clusters. 
%
%
Since each $\text{DMP}_k$ is a linear program (LP), we utilize the LP optimality conditions to reformulate 
model \eqref{eq:Abstract_StabClust} as follows: 
\begin{subequations}\label{eq:StabClust}
\begin{align}
  \underset{\{(\bc^\ell,\cG^\ell)\}_{\ell\in \cL},\{(\bx^{k},\by^k)\}_{k\in \cK}}{\text{minimize}}  &\quad  \underset{\ell\in \cL}\max \quad \underset{k\in \cG^\ell}\max \;  \{d(\hat\bx^k,\bx^{k})\} \label{eq:StabClust_1}\\
 \text{subject to}      & \quad {\bA^k}'\by^{k}=\bc^\ell, \quad \forall \ell \in \cL, k\in \cG^\ell,\label{eq:StabClust_3}\\
	& \quad \by^{k} \geq \bzero, \quad \forall k\in \cK, \label{eq:StabClust_4}\\
	& \quad \bA^k \bx^{k}\geq \bbb^k, \quad \forall k\in \cK,\label{eq:StabClust_5}\\
   	& \quad {\bc^{\ell}}'\bx^{k}={\bbb^k}'\by^k, \quad \quad \forall \ell \in \cL, k\in \cG^\ell, \label{eq:StabClust_6}\\
   	&\quad \lVert\bc^\ell\rVert_1=1, \quad \forall \ell\in \cL. \label{eq:StabClust_2}
\end{align}
\end{subequations}
Constraints \eqref{eq:StabClust_3}--\eqref{eq:StabClust_6} represent the LP optimality conditions for solutions $\{\bx^k\}_{k\in \cG^\ell}$ with respect to $\bc^\ell$ for each cluster $\ell\in \cL$: constraints \eqref{eq:StabClust_3}--\eqref{eq:StabClust_4} enforce dual feasibility where $\by^k\in \mathbb{R}^{m_k}$ represents the vector of dual variables corresponding to $\text{DMP}_k$, constraint \eqref{eq:StabClust_5} corresponds to primal feasibility, and constraint \eqref{eq:StabClust_6} ensures strong duality. Note that problem \eqref{eq:StabClust} is non-convex due to its objective function and constraints \eqref{eq:StabClust_6} and \eqref{eq:StabClust_2}. In Section \ref{sec:sol_struct_reform}, we analyze its solution structure and propose MIP formulations that provide lower and upper bounds on the optimal objective value of problem~\eqref{eq:StabClust}. 
Our subsequent analysis for the rest of this paper focuses on $d(\hat\bx,\bx)=\|\hat \bx-\bx\|_r$ for $r\ge 1$, though similar analysis can be derived for other distance functions. 

\vspace{-0.1in}
\subsection{Heuristics: Two-Stage Approaches}\label{sec:two_stage}
While \eqref{eq:StabClust} provides an exact reformulation of the SC problem, it is computationally challenging. Instead, one naive view on this problem would be to treat the clustering and cost vector inference parts separately. In this subsection, we propose two heuristics based on this separation idea. 

The first algorithm applies traditional K-means clustering to cluster dataset $\hat\cX$ \textit{a priori} based on some distance function, e.g., Euclidean distance, followed by applying inverse optimization post-hoc to derive a cost vector for each of the predetermined clusters. This approach, which we call the cluster-then-inverse (CI) approach, can be written as follows.
%
%
\begin{align}\label{eq:model_CI}
\textbf{CI}(\cK,L):
\begin{cases}
\text{Stage 1}. & \text{Find } \{\cG_{\textup{CI}}^\ell\}_{\ell\in \cL} \in \displaystyle \argmin_{\{\cG^\ell\}_{\ell\in \cL}} \bigg\{\sum_{\ell\in\cL}\sum_{k\in \cG^\ell} d(\bxhat^k,\bx^\ell_\cen)\,\bigg|\,\bx^\ell_\cen \text{ is the centroid of $\cG^\ell$}\bigg\} \\
%
\text{Stage 2}. & \text{Find } \bc^{\ell}_{\textup{CI}} \in \displaystyle \argmin_{\bc^{\ell}}  \bigg\{\underset{k\in \cG_{\textup{CI}}^\ell}{\max}\, d(\hat\bx^k,\bx^k) \,\bigg{|}\, \bx^k\in
\cX^{k*}(\bc^\ell), \lVert\bc^\ell\rVert_1=1\bigg\}, \; \forall \ell\in \cL.
\end{cases}
\end{align}
\noindent In Stage 1, clusters $\{\cG_{\textup{CI}}^\ell\}_{\ell\in \cL}$ are obtained by solving a traditional clustering problem on $\hat\cX$. Then, Stage 2 finds a cost vector for each of the clusters that minimizes cluster instability
. Note that Stage 2 of the CI approach solves a ``reduced'' version of the SC problem that finds a stability-maximizing cost vector for each $\ell\in\cL$ with respect to the observations assigned to cluster $\ell$; i.e., \textbf{SC}$(\cG_{\textup{CI}}^\ell, L=1)$ where $L=1$ implies that no further clustering happens.

Alternatively, the second approach finds a cost vector $\bc^{k*}$ for each data point $k\in \cK$ \textit{a priori} such that 
{$\max \{d(\hat\bx^k,\bx^k)\, |\,  \bx^k\in\cX^{k*}(\bc^{k*})\}$} is minimized. Then, the cost vectors are clustered post-hoc into $L$ groups via traditional clustering
. We call this approach inverse-then-cluster (IC): 
\begin{align}\label{eq:model_IC}
\textbf{IC}(\cK,L):
\begin{cases}
\text{Stage 1.} & \text{Find } \bc^{k*}\in \displaystyle \argmin_{\bc^k} \bigg\{ \max\; d(\hat\bx^k,\bx^k) \,\bigg{|}\, \bx^k\in\cX^{k*}(\bc^k),\,  \lVert\bc^k\rVert_1=1 \bigg\}, \; \forall k\in \cK\\
\text{Stage 2.} & \text{Find } \{\cG_{\textup{IC}}^\ell\}_{\ell\in \cL} \in \displaystyle \argmin_{\{\cG^\ell\}_{\ell\in \cL}} \bigg\{  \sum_{\ell\in\cL}\sum_{k\in \cG^\ell}d(\bc^{k*}, \bc^\ell_\cen) \,\bigg{|}\, \bc^\ell_\cen \text{ is the centroid of } \cG^\ell \bigg\}. %
\end{cases}
\end{align}
Note that, similarly, Stage 1 of the above IC approach can be seen as solving a reduced version of SC, i.e., \textbf{SC}$(\{k\}, L=1)$, which finds a ``per-observation'' cost vector $\bc^{k*}$ that maximizes stability associated with each observation $\hat\bx^k$. However, once the cost vectors are clustered in Stage 2, it is the resulting centroid cost vector, i.e., $\bc^\ell_\cen$, that represents the preferences for the observations assigned to cluster $\ell$, which does not necessarily retain the same level of stability achieved by the per-observation cost vectors (i.e., $\bc^{k*}$'s) in Stage 1. To address this, once the clustering is done, one may solve the SC problem for each cluster again to find a ``corrected'' cost vector; i.e., \textbf{SC}$(\cG_{\textup{IC}}^\ell,L=1)$ for each $\ell \in \cL$. We denote such a post-processed cost vector by $\bc^{\ell}_{\textup{IC}},\ell \in \cL$.


\subsection{Model Comparison}
Next, we compare the performance of the SC model (i.e., \eqref{eq:Abstract_StabClust}) and the CI and IC approaches.
\begin{proposition} \label{prop:models_compare}
Given $\hat\cX,$ let $\{\cG_{\textup{SC}}^{\ell},\bc^{\ell}_{\textup{SC}}\}_{\ell\in \cL}$ denote an optimal solution to model \eqref{eq:Abstract_StabClust}, and $\{\cG_{\textup{CI}}^{\ell},\bc^{\ell}_{\textup{CI}}\}_{\ell\in \cL}$ and $\{\cG_{\textup{IC}}^{\ell},\bc^{\ell}_{\textup{IC}}\}_{\ell\in \cL}$ be the clusters and corresponding cost vectors achieved by the CI and IC approaches, respectively. Then we have
\begin{itemize}
\item [(i)] $\underset{\ell\in \cL,k\in \cG_{\textup{SC}}^\ell}{\max}\{ d(\hat\bx^k,\bx^k)\ |\ \bx^k\in\cX^{k*}(\bc^\ell_{\textup{SC}})\} \le \underset{\ell\in \cL,k\in \cG_{\textup{CI}}^\ell}{\max}\{d(\hat\bx^k,\bx^k)\ |\ \bx^k\in\cX^{k*}(\bc^\ell_{\textup{CI}})\} $, and 
\item [(ii)]$\underset{\ell\in \cL,k\in \cG_{\textup{SC}}^\ell}{\max}\{d(\hat\bx^k,\bx^k)\ |\ \bx^k\in\cX^{k*}(\bc^\ell_{\textup{SC}})\} \le \underset{\ell\in \cL,k\in \cG_{\textup{IC}}^\ell}{\max}\{d(\hat\bx^k,\bx^k)\ |\ \bx^k\in\cX^{k*}(\bc^\ell_{\textup{IC}})\}$.
\end{itemize}
\end{proposition}

\begin{proofin} Since $\{\cG_{\textup{SC}}^{\ell},\bc^{\ell}_{\textup{SC}}\}_{\ell\in \cL}$ is an optimal solution to \eqref{eq:Abstract_StabClust}, we have\begin{align}\nonumber
\underset{\ell\in \cL,k\in \cG_{\textup{SC}}^\ell}{\max}\{ d(\hat\bx^k,\bx^k)\ |\ \bx^k\in\cX^{k*}(\bc^\ell_{\textup{SC}})\}&=\underset{\{(\bc^\ell,\cG^\ell)\}_{\ell\in \cL}}{\min} \bigg\{\underset{\ell\in \cL,k\in \cG^\ell}{\max}\; \{ d(\hat\bx^k,\bx^k)\ |\ \bx^k\in\cX^{k*}(\bc^\ell)\} \,\bigg|\, \|\bc^\ell\|_1=1\bigg\},
\end{align}
where the right hand side corresponds to model \eqref{eq:Abstract_StabClust}. 

For part (i), consider $\{\cG^\ell_{\textup{CI}}, \bc^{\ell}_{\textup{CI}}\}_{\ell\in \cL}$ generated by CI. 
Recall from \eqref{eq:model_CI} that $\bc^{\ell}_{\textup{CI}} \in \displaystyle \argmin_{\bc^{\ell}} \bigg\{\underset{k\in \cG_{\textup{CI}}^\ell}{\max}\;  \{d(\hat\bx^k,\bx^k)\, | \, \bx^k\in\cX^{k*}(\bc^\ell)\}   \,\bigg{|}\,  \lVert\bc^\ell\rVert_1=1\bigg\}$ for the given cluster $\cG^\ell_{\textup{CI}}$ for each $\ell$. Thus, $\underset{\ell\in \cL,k\in \cG_{\textup{CI}}^\ell}{\max}\{d(\hat\bx^k,\bx^k)\, |\, \bx^k\in\cX^{k*}(\bc^\ell_{\textup{CI}})\} = \underset{\ell\in\cL}{\max}\; \underset{\bc^{\ell}}{\min} \bigg\{\underset{k\in \cG_{\textup{CI}}^\ell}{\max} \{d(\hat\bx^k,\bx^k)\, | \, \bx^k\in\cX^{k*}(\bc^\ell)\}   \,\bigg{|}\,  \lVert\bc^\ell\rVert_1=1\bigg\}.$ Then it follows that 
\begin{align}\nonumber
& \underset{\ell\in \cL,k\in \cG^\ell_{\textup{SC}}}{\max}\{ d(\hat\bx^k,\bx^k)\, |\, \bx^k\in\cX^{k*}(\bc^\ell_{\textup{SC}})\} \nonumber\\
& = \underset{\{(\bc^\ell,\cG^\ell)\}_{\ell\in \cL}}{\min} \bigg\{\underset{\ell\in \cL,k\in \cG^\ell}{\max}\; \{ d(\hat\bx^k,\bx^k)\, |\, \bx^k\in\cX^{k*}(\bc^\ell)\} \,\bigg|\, \|\bc^\ell\|_1=1\bigg\} \nonumber\\
& \le \underset{\ell\in\cL}{\max}\; \underset{\bc^{\ell}}{\min} \bigg\{\underset{k\in \cG_{\textup{CI}}^\ell}{\max}\;  \{d(\hat\bx^k,\bx^k)\, | \, \bx^k\in\cX^{k*}(\bc^\ell)\}   \,\bigg{|}\,  \lVert\bc^\ell\rVert_1=1\bigg\} \nonumber\\
& = \underset{\ell\in \cL,k\in \cG^\ell_{\textup{CI}}}{\max}\; \{ d(\hat\bx^k,\bx^k)\, |\, \bx^k\in\cX^{k*}(\bc^\ell_{\textup{CI}})\},\nonumber
\end{align}
as desired.

Proof for part (ii) is similar.  Consider $\{\cG^\ell_{\textup{IC}},\bc^{\ell}_{\textup{IC}}\}_{\ell\in \cL}$ generated by IC. 
Recall from \eqref{eq:model_IC} and its post-processing step that $\bc^{\ell}_{\textup{IC}} \in \displaystyle \argmin_{\bc^{\ell}} \bigg\{\underset{k\in \cG_{\textup{IC}}^\ell}{\max}\;  \{d(\hat\bx^k,\bx^k)\, |\, \bx^k\in\cX^{k*}(\bc^\ell)\} \,  \bigg{|}  \, \lVert\bc^\ell\rVert_1=1\bigg\}$. 
That is, 
$\underset{\ell\in \cL,k\in \cG^\ell_{\textup{IC}}}{\max}\!\!\{ d(\hat\bx^k,\bx^k)\, |\, \bx^k\in\cX^{k*}(\bc^\ell_{\textup{IC}})\} = \underset{\ell\in\cL}{\max}\,\underset{\bc^{\ell}}{\min} \bigg\{\underset{k\in \cG_{\textup{IC}}^\ell}{\max} \{d(\hat\bx^k,\bx^k)\, |\, \bx^k\in\cX^{k*}(\bc^\ell)\} \,  \bigg{|}  \, \lVert\bc^\ell\rVert_1=1\bigg\}$. 
Thus, we have 
\begin{align}
&\underset{\ell\in \cL,k\in \cG^\ell_{\textup{SC}}}{\max} \{ d(\hat\bx^k,\bx^k)\, |\, \bx^k\in\cX^{k*}(\bc^\ell_{\textup{SC}})\} \nonumber\\
&= \underset{\{(\bc^\ell,\cG^\ell)\}_{\ell\in \cL}}{\min} \bigg\{\underset{\ell\in \cL,k\in \cG^\ell}{\max} \{ d(\hat\bx^k,\bx^k)\, |\, \bx^k\in\cX^{k*}(\bc^\ell)\} \,\bigg|\, \|\bc^\ell\|_1=1\bigg\} \nonumber\\
&\le \underset{\ell\in\cL}{\max}\,\underset{\bc^{\ell}}{\min} \bigg\{\underset{k\in \cG_{\textup{IC}}^\ell}{\max}\;  \{d(\hat\bx^k,\bx^k)\, |\, \bx^k\in\cX^{k*}(\bc^\ell)\} \,  \bigg{|}  \, \lVert\bc^\ell\rVert_1=1\bigg\} \nonumber \\
&=\underset{\ell\in \cL,k\in \cG^\ell_{\textup{IC}}}{\max} \{ d(\hat\bx^k,\bx^k)\, |\, \bx^k\in\cX^{k*}(\bc^\ell_{\textup{IC}})\},\nonumber
\end{align}
as desired. \hfill$\square$
\end{proofin}

While Proposition~\ref{prop:models_compare} implies that the SC model performs at least as well as CI and IC in terms of stability, the SC model is typically computationally more challenging than CI and IC. 
In the next section, we analyze the solution structure of the SC model, which we use to derive MIP formulations that provide lower and upper bounds on the optimal value of the SC model. %


\vspace{-0.1in}
\section{Solution Structure and Bounds}\label{sec:sol_struct_reform}
\vspace{-0.1in}
The reformulation of the SC model (i.e., \eqref{eq:StabClust}) is non-convex due to the normalization constraint \eqref{eq:StabClust_2} as well as the objective function: for a given $k\in \cK$ and arbitrary $\bc$, 
$\max  \{d(\hat\bx^k,\bx^k)\, |\, \bx^k\in \cX^{k*}(\bc)\}$ is a maximization of the convex function $d$ over the convex region $\cX^{k*}(\bc)$. Both the CI and IC approaches also face the same computational challenges because they also involve solving the SC formulations albeit of smaller size; i.e., Stage 2 of the CI approach solves $\textbf{SC}(\cG_{\text{CI}}^\ell,L=1)$ for each 
$\ell\in \cL$ 
and Stage 1 of the IC approach solves $\textbf{SC}(\{k\},L=1)$ for each $k\in\cK$. 
In this section, we analyze the solution structure of the SC model, which leads to MIP formulations 
that provide lower and upper bound solutions for the SC problem.
\begin{theorem} \label{theorem:corner_point}
There exists an optimal solution $\big(\{(\bc^{\ell*},\cG^{\ell*})\}_{\ell\in \cL},\{(\bx^{k*},\by^{k*})\}_{k\in \cK}\big)$ to \eqref{eq:StabClust} such that for each cluster $\ell\in
\cL$:  
\begin{itemize}
    \item [\textup{(i)}] ${\ba^{ki}}'\bx^{k*}=b_i^k$ for $i\in \cI^{k*}\subseteq \cI^k$ where $|\cI^{k*}|=n$ for all $k\in \cG^{\ell*}$, and 
    \item [\textup{(ii)}] $\bc^{\ell*}\in \textup{cone}(\{\ba^{ki}\}_{i\in \cI^{k*}})$ for all $k\in \cG^\ell$ where $\textup{cone}(\cdot)$ denotes the conic hull of the given vectors, i.e, $\textup{cone}(\{\ba^{ki}\}_{i\in \cI^{k*}})=\{\sum_{i\in \cI^{k*}} \gamma_i\ba^{ki}\ | \ \gamma_i\ge 0\}$.
\end{itemize}
\end{theorem}
\begin{proofin} 
Consider an optimal solution $\big(\{(\bc^{\ell*},\cG^{\ell*})\}_{\ell\in \cL},\{(\bx^{k*},\by^{k*})\}_{k\in \cK}\big)$ to \eqref{eq:StabClust}.  
Due to constraints \eqref{eq:StabClust_3}--\eqref{eq:StabClust_6}, we have $\bx^{k*}\in \cX^{k*}(\bc^{\ell*})
$ for each $\ell\in \cL$ and $k\in \cG^{\ell*}$. Note that any point in $\cX^{k*}(\bc^{\ell*})$ can be represented by a convex combination of extreme points of $\cX^{k*}(\bc^{\ell*})$. Let $\text{ext}(\cX^{k*}(\bc^{\ell*}))$ be the set of extreme points of $\cX^{k*}(\bc^{\ell*})$, $Q_k=\left|\text{ext}(\cX^{k*}(\bc^{\ell*}))\right|$, and $\cQ^k=\{1,\ldots,Q_k\}$, i.e., $\text{ext}(\cX^{k*}(\bc^{\ell*}))=\{\bar\bx^1,\ldots, \bar\bx^{Q_k}\}$, for each $k\in \cK$. 
Then, there exists $\bar\blambda\in \mathbb{R}_+^{Q_k}$ such that $\bx^{k*}=\sum_{q_k\in \cQ^k} \bar\lambda_{q_k}\bar\bx^{q_k}$ and $\sum_{q_k\in \cQ^k}\bar\lambda_{q_k}=1$. 

Now we prove part (i). Let $q^*_k\in \displaystyle \argmax_{q_k\in \cQ^k} \{\|\hat\bx^k-\bar\bx^{q_k}\|_r\}$. That is, we have $\|\hat\bx^k-\bar\bx^{q^*_k}\|_r\ge \|\hat\bx^k-\bar\bx^{q_k}\|_r$ for all $q_k\in \cQ^k$. Multiplying both sides of the inequality by $\bar\lambda_{q_k}$ yields $\bar\lambda_{q_k}\|\hat\bx^k-\bx^{q^*_k}\|_r\ge \bar\lambda_{q_k}\|\hat\bx^k-\bar\bx^{q_k}\|_r$ for all $q_k\in \cQ^k$, and thus $\displaystyle\sum_{q_k\in \cQ^k}\bar\lambda_{q_k}\|\hat\bx^k-\bar\bx^{q^*_k}\|_r \ge \displaystyle\sum_{q_k\in \cQ^k}\bar\lambda_{q_k}\|\hat\bx^k-\bar\bx^{q_k}\|_r$. Note that, from $\displaystyle\sum_{q_k\in \cQ^k}\bar\lambda_{q_k}= 1$ we have $\displaystyle\sum_{q_k\in \cQ^k}\bar\lambda_{q_k}\|\hat\bx^k-\bar\bx^{q^*_k}\|_r=\|\hat\bx^k-\bar\bx^{q^*_k}\|_r$. This leads to 
\begin{align}\nonumber
\|\hat\bx^k-\bar\bx^{q^*_k}\|_r 
& \ge \displaystyle\sum_{q_k\in \cQ^k}\bar\lambda_{q_k}\|\hat\bx^k-\bar\bx^{q_k}\|_r \nonumber\\
&= \displaystyle\sum_{q_k\in \cQ^k}\|\bar\lambda_{q_k} \hat\bx^k-\bar\lambda_{q_k}\bar\bx^{q_k}\|_r\nonumber \\
&\ge \|\displaystyle\sum_{q_k\in \cQ^k}\bar\lambda_{q_k}\hat\bx^k-\displaystyle\sum_{q_k\in \cQ^k}\bar\lambda_{q_k}\bar\bx^{q_k}\|_r \nonumber\\
&= \|\hat\bx^k-\displaystyle\sum_{q_k\in \cQ^k}\bar\lambda_{q_k}\bar\bx^{q_k}\|_r \nonumber\\
&= \|\hat\bx^k-\bx^{k*}\|_r, \nonumber
\end{align}
where the second inequality holds due to Minkowski inequalities. 
%
Also, from the optimality of $\bx^{k*}$, we have $\|\hat\bx^k-\bar\bx^{q^*_k}\|_r \le \|\hat\bx^k-\bx^{k*}\|_r$
. Thus, it must be that $\|\hat\bx^k-\bar\bx^{q^*_k}\|_r= \|\hat\bx^k-\bx^{k*}\|_r$. This means that the solution $\big(\{(\bc^{\ell*},\cG^{\ell*})\}_{\ell\in \cL},\{(\bar\bx^{q^*_k},\by^{k*})\}_{k\in \cK}\big)$ is also optimal to \eqref{eq:StabClust}. Since $\bar\bx^{q^*_k}$ is an extreme point, there must exist $\cI^{k*}\subseteq\cI^k$ such that ${\ba^{ki}}'\bar\bx^{q_k^*}=b_i^k$ for all $i\in \cI^{k*}$ and $|\cI^{k*}|=n$. 

We prove part (ii) using 
 the same above optimal solution $\big(\{(\bc^{\ell*},\cG^{\ell*})\}_{\ell\in \cL},\{(\bar\bx^{q^*_k},\by^{k*})\}_{k\in \cK}\big)$ to \eqref{eq:StabClust}
. First, note that for each $\ell\in \cL$, $\bc^{\ell*}$ satisfies \eqref{eq:StabClust_2}, which means for all $k\in \cG^{\ell*}$ there exists at least one $i\in\cI^{k*}$ for which $y^{k*}_i>0$. 
Moreover, because ${\ba^{ki}}'\bar\bx^{q^*_k}> b_i$ for $i\in \cI^k\setminus \cI^{k*}$ and $y^{k*}_i$ is the associated dual variable, it must be that $y^{k*}_i=0$ for all $i\in \cI^k\setminus \cI^{k*}$. Thus, from \eqref{eq:StabClust_3} we have $\displaystyle\bc^{\ell*}=\sum_{i\in \cI^k}y^{k*}_i \ba^{ki}=\sum_{i\in \cI^{k*}}y^{k*}_i \ba^{ki}$, or equivalently $\bc^{\ell*}\in \text{cone}(\{\ba^{ki}\}_{i\in \cI^{k*}})$, for all $k\in \cG^{\ell*}$. \hfill$\square$
\end{proofin}

The following result characterizes the solution structure of the SC model under the special case where all DMs solve the same DMP. 
\begin{corollary}\label{corollary:same_As}
Assume $\bA^k=\bA$ and $\bbb^k=\bbb$ for all $k\in \cK$ and let $\cI$ be the index set for rows of $\bA$. Then there exists an optimal solution $\big(\{(\bc^{\ell*},\cG^{\ell*})\}_{\ell\in \cL},\{(\bx^{k*},\by^{k*})\}_{k\in \cK}\big)$ to \eqref{eq:StabClust} such that $\bc^{\ell*}\in \textup{cone}_+(\{\ba^{i}\}_{i\in \cI^*})$ for each cluster $\ell\in
\cL$, where $\cI^*\subseteq\cI$, $|\cI^*|=n$, and $\textup{cone}_+(\cdot)$ denotes the interior of the conic hull of given vectors, i.e., $\textup{cone}_+(\{\ba^{i}\}_{i\in \cI^*})=\{\sum_{i\in \cI^{*}} \lambda_i\ba^{i}\, | \, \lambda_i>0,\forall i\in\cI^*\}$.
\end{corollary}

\vspace{-0.1in}
\subsection{Lower Bound Formulation}
\vspace{-0.1in}Theorem \ref{theorem:corner_point} states that there exists an optimal solution to the SC model where $\bx^{k*}$ is an extreme point of $\cX^k$ for all $k\in \cK$. Also, if $k\in \cG^{\ell*}$ then $\bc^{\ell*}$ must be a conic combination of $\ba^{ki}$'s for $i$ such that ${\ba^{ki}}'\bx^{k*}=b_i^k$. Based on this observation, we propose an MIP formulation that explicitly finds an extreme point $\bx^{k*}$ for each $\cX^{k}$, clusters the data points, and constructs $\bc^{\ell*}$ for cluster $\ell$ as a conic combination of $\ba^{ki}$'s for $k$ assigned to cluster $\ell$ and for $i$ such that ${\ba^{ki}}'\bx^{k*}=b_i^k$. We then show that the optimal value of this MIP is a lower bound on that of the SC problem:
\begin{subequations}\label{eq:Re_StabClust_MIP}
\begin{align}
\textbf{SC-LB}(\cK,\cL):  
  \underset{\underset{\{(\bc^{\ell},\bc^{\ell+},\bc^{\ell-},\bz^\ell)\}_{\ell\in \cL},\bu}{\{(\bx^k,\bv^{k},\blambda^{k})\}_{k\in \cK},\textcolor{white}{\alpha}}}{\text{minimize}}  &  \quad \underset{k\in\cK}{\max}\{d(\hat\bx^k,\bx^{k})\}\\
 \text{subject to}           
     & \quad {\bA^k}'\blambda^{k} -M_1(1-u_{k\ell})\le \bc^\ell \nonumber\\
     & \hspace{0.3in}\le {\bA^k}'\blambda^{k}+M_1(1-u_{k\ell}), \quad \forall \ell \in \cL,k\in \cK,  \label{eq:Re_StabClust_MIP_6}\\
     & \quad \lambda^k_i\le M_2 v^k_i, \quad \forall k\in \cK, i\in \cI^k,   \label{eq:Re_StabClust_MIP_7}\\
     & \quad b^k_i\le {\ba^{ki}}'\bx^{k} \leq b_i^k +M_3(1-v^k_i),\quad \forall k \in \cK, i  \in \cI^k, \label{eq:Re_StabClust_MIP_8}\\
     & \quad \sum_{i\in \cI^k} v^k_i=n, \quad \forall k\in \cK,  \label{eq:Re_StabClust_MIP_9}\\
     & \quad \sum_{\ell\in \cL} u_{k\ell}=1, \quad  \forall k\in \cK,   \label{eq:Re_StabClust_MIP_10}\\
          &\quad \bc^\ell=\bc^{\ell+}-\bc^{\ell-}, \quad \forall \ell\in \cL,\label{eq:Re_StabClust_MIP_2}\\
     &\quad \bc^{\ell+}\le \bz^\ell, \quad \forall \ell\in \cL,\label{eq:Re_StabClust_MIP_3}\\
     &\quad \bc^{\ell-}\le \beee-\bz^\ell, \quad \forall \ell\in \cL,\label{eq:Re_StabClust_MIP_4}\\
     &\quad \beee'(\bc^{\ell+}+\bc^{\ell-})=1,\quad \forall \ell\in \cL,\label{eq:Re_StabClust_MIP_5}\\
     &\quad \bv^k\in \{0,1\}^{n}, \bu\in \{0,1\}^{K\times L}, \bz^\ell\in \{0,1\}^n, \quad  \forall k\in \cK,\ell \in \cL,  \label{eq:Re_StabClust_MIP_12}\\
          &\quad \blambda^k, \bc^{\ell+},\bc^{\ell-}\ge \bzero, \quad \forall k\in \cK, \forall \ell \in \cL,\label{eq:Re_StabClust_MIP_13}
\end{align}
\end{subequations}
where parameters $M_1$, $M_2$, and $M_3$ are sufficiently large positive constants. 
Using the result of Theorem \ref{theorem:corner_point}, constraints \eqref{eq:Re_StabClust_MIP_6}--\eqref{eq:Re_StabClust_MIP_7} enforce each $\bc^\ell$ to be a conic combination of some $\ba^{ki}$'s; which $\ba^{ki}$ is selected is dictated by binary variables $v_i^k$ and $u_{k\ell}$. If $u_{k\ell}=1$, data point $\hat\bx^k$ is assigned to cluster $\ell$ and \eqref{eq:Re_StabClust_MIP_6} holds with equality. 
The variables $\lambda_i^k$ in \eqref{eq:Re_StabClust_MIP_6} are then controlled by \eqref{eq:Re_StabClust_MIP_7} using binary variable $v_i^k$, i.e., if $v_i^k=0$ then $\lambda_i^k=0$ and thus preventing $\ba^{ki}$ from being a basis vector for the conic hull constructing $\bc^\ell$. Constraints \eqref{eq:Re_StabClust_MIP_8}--\eqref{eq:Re_StabClust_MIP_9} enforce each $\bx^k$ to be an extreme point of $\cX^k$, i.e.,  
satisfying ${\ba^{ki}}'\bx^k\ge b^k_i$ with equality for $n$ number of $i$'s ensured by \eqref{eq:Re_StabClust_MIP_9}. Constraint \eqref{eq:Re_StabClust_MIP_10} ensures that each observation is assigned to only one cluster
. Finally, constraints \eqref{eq:Re_StabClust_MIP_2}--\eqref{eq:Re_StabClust_MIP_5} replace the non-convex normalization constraint \eqref{eq:StabClust_2}. 
%
%
%
%
%
%
%
The following result shows that the optimal value of the above problem is a lower bound on the optimal value of the SC problem, i.e., \eqref{eq:StabClust}. 
\begin{proposition} \label{prop:strict_MIP}
Let $\rho^*$ and $ \beta^*$ denote the optimal objective values of problems \eqref{eq:StabClust} and \eqref{eq:Re_StabClust_MIP}, respectively. 
Then, we have
\textup{(i)} $ \beta^*\le \rho^*$, and \textup{(ii)} $\beta^* = \rho^*$ if 
there exists $(\{\tilde\bv^{k}\}_{k\in \cK}, \tilde\bu, \{\tilde\bc^{\ell}\}_{\ell\in \cL})$ optimal for 
\eqref{eq:Re_StabClust_MIP} such that $\tilde\bc^{\ell}\in \textup{cone}_+(\{\ba^{ki}\}_{i,k:\tilde{v}^{k}_i=1, \tilde{u}_{k\ell}=1})$ for each $\ell\in \cL$.
\end{proposition}

Proposition \ref{prop:strict_MIP} suggests that the optimal value of model \eqref{eq:Re_StabClust_MIP} is a lower bound on that of the SC model. Proposition \ref{prop:strict_MIP}  also implies that once model \eqref{eq:Re_StabClust_MIP} is solved, we can check the condition in Proposition \ref{prop:strict_MIP}(ii) to determine whether \eqref{eq:Re_StabClust_MIP} achieves the exact optimal value of the SC model.

\subsection{Upper Bound Formulation}
Proposition \ref{prop:strict_MIP} states that if formulation \eqref{eq:Re_StabClust_MIP} finds a solution such that each 
$\bc^\ell$, $\ell\in \cL$, is a strict conic combination of the selected $\ba^{ki}$ vectors, then its optimal value is equal to that of the SC model. Based on this observation, we add a constraint to \eqref{eq:Re_StabClust_MIP} that enforces this condition and show that the following modified problem provides an upper bound on the optimal value of the SC model:
\begin{subequations}\label{eq:Re_StabClust_MIP_UB}
\begin{align}
\textbf{SC-UB}(\cK,\cL):  \underset{\underset{\{(\bc^{\ell},\bc^{\ell+},\bc^{\ell-},\bz^\ell)\}_{\ell\in \cL},\bu}{\{(\bx^k,\bv^{k},\blambda^{k})\}_{k\in \cK},\textcolor{white}{\alpha}}}{\text{minimize}}  &  \quad \underset{k\in\cK}{\max}\{d(\hat\bx^k,\bx^{k})\}   
  \label{eq:Re_StabClust_MIP_UB_1}\\ 
 \text{subject to} 
     & \quad \eqref{eq:Re_StabClust_MIP_6}-\eqref{eq:Re_StabClust_MIP_13} \label{eq:Re_StabClust_MIP_UB_2}\\
     & \quad \textcolor{black}{\lambda^k_i\ge v^k_i\hat\alpha, \quad \forall k\in \cK, i\in \cI^k,} \label{eq:Re_StabClust_MIP_UB_3}
\end{align}
\end{subequations}
where $\hat \alpha$ is a small positive constant. For each $\ell\in \cL$, $k\in \cG^\ell$, and $i\in \cI^k$, if $v_i^k=1$ then $\lambda_i^k \ge \hat\alpha>0$, which enforces $\bc^\ell$ to be a strict conic combination of $n$ selected $\ba^{ki}$ vectors (i.e., for which $v_i^{k}=1$; see \eqref{eq:Re_StabClust_MIP_8}--\eqref{eq:Re_StabClust_MIP_9}). If there exists an optimal solution to the SC problem whose $\bc^\ell$ vectors satisfy the strict conic combination condition then \eqref{eq:Re_StabClust_MIP_UB} with an appropriate $\hat\alpha$ generates the optimal solution for the SC model; otherwise, the optimal value of \eqref{eq:Re_StabClust_MIP_UB} is an upper bound on the optimal value of the SC model. We formalize this in the following result.
\begin{proposition} \label{prop:strict_MIP_UB}
Given $\hat\alpha$, let $\rho^*$ and $ \beta^*$ denote the optimal objective values of problems \eqref{eq:StabClust} and \eqref{eq:Re_StabClust_MIP_UB}, respectively. Then, we have $\rho^*\le \beta^*$.
\end{proposition}

While \textbf{SC-LB}($\cK,\cL$) and \textbf{SC-UB}($\cK,\cL$) provide bounds for the SC problem, i.e., \textbf{SC}($\cK,\cL$), these problems are typically large-scale MIPs and thus can be computationally challenging. Appendix~\ref{sec:solution_approaches} shows how the CI and IC approaches can be used to create  initial feasible solutions for these MIPs and reduce the computational burden.


\vspace{-0.1in}
\section{Numerical Results}\label{sec:Results}
\vspace{-0.1in}
\rev{In this section, we first examine the performance of the proposed clustering approach using various-sized instances and discuss the computational benefits and limitations. We then present the results of the application of the proposed approach in the diet recommendation context to cluster DMs based on the similarity of their food preferences.}

\subsection{\rev{Performance of the Proposed Clustering Approaches}}\label{sec:results_synthetic}
We use various-sized randomly generated instances to demonstrate the CI, IC, and SC approaches. For small instances we chose $K\in \{30,40,50\}$ and generated LP instances with $n=10$ and $m^{k} = m\in \{30,40\}, \forall k=1,\ldots,K$. For large instances we used $K\in \{100,115,130\}$, $n=20$, and $m^{k} =  m\in \{60,80\}, \forall k=1,\ldots,K$. To generate dataset $\hat\cX$ for each instance, we generated $K$ random cost vectors, solved $K$ DMPs to generate optimal solutions, and added random noise to the solutions. All optimization problems were solved by Gurobi 9.1 \citep{gurobi} with a 16-core 2.9 GHz processor and 512 GB memory.

\begin{table}[htbp]
  \centering
  \caption{Performance of IC, CI, and SC approximated by upper and lower bounds.}
    \label{tab:Results}%
    \begin{tabular}{ccrrrrrrrrr}
    \toprule
          &       & \multicolumn{4}{c}{Worst-case distance}   &       & \multicolumn{4}{c}{Time (s)} \\
        \cmidrule{3-6}\cmidrule{8-11}
\multicolumn{1}{l}{$(n,m)$} & $K$     & \multicolumn{1}{c}{IC} & \multicolumn{1}{c}{CI} & \multicolumn{1}{c}{UB} & \multicolumn{1}{c}{LB} &       & \multicolumn{1}{c}{IC} & \multicolumn{1}{c}{CI} & \multicolumn{1}{c}{UB} & \multicolumn{1}{c}{LB} \\
    \midrule
    \multirow{3}[2]{*}{(10,30)} & 30    & 14.71 & 11.67 & 1.92  & 1.92  &       & 36.48 & 39.01 & 70.77 & 13.36 \\
          & 40    & 13.09 & 13.97 & 1.78  & 1.77  &       & 112.55 & 174.39 & 273.01 & 20.31 \\
          & 50    & 14.58 & 14.26 & 2.01  & 1.86  &       & 147.84 & 216.58 & 5707.04 & 27.75 \\
    \midrule
    \multirow{3}[2]{*}{(10,40)} & 30    & 9.51  & 9.57  & 1.97  & 1.97  &       & 74.84 & 94.52 & 238.45 & 16.43 \\
          & 40    & 10.81 & 11.17 & 1.50  & 1.50  &       & 158.76 & 119.01 & 470.01 & 26.47 \\
          & 50    & 12.27 & 12.15 & 2.09  & 1.99  &       & 575.92 & 321.21 & 1667.48 & 40.78 \\
    \midrule
    \multirow{3}[2]{*}{(20,60)} & 100   & 17.04 & 15.40 & 2.03  & 1.97  &       & 675.62 & 443.55 & 4441.23 & 426.37 \\
          & 115   & 16.13 & 19.40 & 2.08  &   1.98    &       & 789.73 & 517.47 & 3417.01 & 589.38 \\
          & 130   & 17.86 & 16.27 & 2.01  & 1.97  &       & 917.40 & 549.12 & 4323.03 & 691.85 \\
    \midrule
    \multirow{3}[2]{*}{(20,80)} & 100   & 15.20 & 13.31 & 2.03  & 1.85  &       & 833.62 & 609.72 & 2874.12 & 807.51 \\
          & 115   & 14.98 & 17.41 & 2.11  &  1.92    &       & 885.66 & 739.27 & 4404.16 & 991.61 \\
          & 130   & 14.50 & 13.85 & 2.06  & 2.00  &       & 1667.74 & 1165.07 & 5720.49 & 1141.38 \\
    \bottomrule
    \end{tabular}%
\end{table}

Table \ref{tab:Results} shows the worst-case distances and solution times achieved by the IC and CI approaches as well as the upper and lower bound formulations for the SC problem (i.e., \textbf{SC-UB}($\cK,\cL$) and \textbf{SC-LB}($\cK,\cL$), respectively). Although the IC and CI problems involved solving smaller versions the SC problem, which were approximated by solving their respective smaller versions of both SC-UB and SC-LB problems, for brevity Table \ref{tab:Results} only presents the IC and CI results approximated by the smaller version of SC-UB. Columns labeled UB in Table \ref{tab:Results}  show the results for \textbf{SC-UB}($\cK,\cL$), which were obtained using an initial solution achieved by the IC and CI results presented in this table (see Appendix~\ref{sec:solution_approaches}); thus, the solution time for UB is the time for finding an initial solution via either IC or CI (whichever gives a smaller worst-case distance) plus the time for the solver to improve the initial solution and find an optimal solution. Columns labeled LB show the results for \textbf{SC-LB}($\cK,\cL$). For each instance $(n,m,K)$, the reported worst-case distance values and times in the table were averaged over two sub-instances with $L=3$ and $L=5$. For all instances, we can see that the UB and LB values were close to each other or identical, indicating that the solutions from both the upper and lower bound formulations are close to the optimal solutions to the SC model. Since the SC model considers a minimization of the worst-case distance, our suggestion is to use the clusters and cost vectors achieved by the upper bound formulation so as not to underestimate the true cluster instability. 

\begin{table}[htbp]
  \centering
  \caption{Performance of IC, CI, and SC for instances with the same $(\bA,\bbb)$ for all DMs.}
    \begin{tabular}{ccrrrrrrr}
\toprule         &       & \multicolumn{3}{c}{Worst-case distance} &       & \multicolumn{3}{c}{Time (s)} \\
\cmidrule{3-5}\cmidrule{7-9}    \multicolumn{1}{l}{$(n,m)$} & $K$     & \multicolumn{1}{c}{IC} & \multicolumn{1}{c}{CI} & \multicolumn{1}{c}{SC} &       & \multicolumn{1}{c}{IC} & \multicolumn{1}{c}{CI} & \multicolumn{1}{c}{SC} \\
    \midrule
    \multirow{3}[2]{*}{(5,15)} & 20    & 2.21  & 1.96  & 1.96  &       & 4.96  & 3.15  & 4.45 \\
          & 30    & 3.48  & 1.99  & 1.99  &       & 9.69  & 6.19  & 8.56 \\
          & 40    & 2.91  & 1.95  & 1.95  &       & 15.64 & 11.71 & 15.20 \\
    \midrule
    \multirow{3}[2]{*}{(10,30)} & 20    & 2.03  & 2.17  & 1.99  &       & 9.31  & 6.38  & 13.73 \\
          & 30    & 3.93  & 1.96  & 1.96  &       & 72.92 & 16.93 & 69.88 \\
          & 40    & 3.91  & 2.09  & 2.02  &       & 99.84 & 32.24 & 96.15 \\
    \bottomrule
    \end{tabular}%
  \label{tab:special_instances}%
\end{table}%

The performance of the CI and IC approaches depends highly on the geometric variation of the DMP feasible regions. For example, when all DMs solve DMPs with similar constraints, i.e., similar $\bA^k$ and $\bbb^k$, the performance of CI and IC becomes comparable to that of the SC approach. To demonstrate this, we generated instances where $(\bA,\bbb)$ is fixed across all DMs. Table \ref{tab:special_instances} shows the result of CI, IC, and SC for these instances. Recall from Corollary \ref{corollary:same_As} that solving the upper bound formulation \textbf{SC-UB}($\cK,\cL$), i.e., \eqref{eq:Re_StabClust_MIP_UB}, for these instances generates an optimal solution for the SC model. In fact, we solved both \eqref{eq:Re_StabClust_MIP_UB} (with $\hat\alpha=0.05$) and \eqref{eq:Re_StabClust_MIP} and their optimal values matched, indicating that the solution is indeed optimal for SC. In most cases, both IC and CI find the worst-case distance values close to those from the SC model, though CI appears to perform better than IC for these specific instances. 
\subsection{\rev{Application to the Diet Problem: Clustering Based on Food Preferences}}\label{sec:results_diet}

\rev{
Clustering has been widely used in the diet and nutrition literature for identifying distinct diet patterns from a specific subject group, associating them with individual health conditions and socio-demographic indicators, and predicting future diets or recommending alternative healthy diets \cite{baek2019hybrid, brennan2010dietary, mcnaughton2008dietary, newby2004empirically}. The existing clustering approach focuses on the similarity of food choices themselves within a homogeneous subject group (e.g., children, diabetic patients, etc.); this, however, often fails to capture unique preferences of the individuals (or DMs) if there is any variation in the underlying nutritional or budgetary requirements among different DMs; for example, the exact same diet patterns might be viewed very differently under different nutritional requirements between diabetic patients and others. 
Recently, inverse optimization has been used to leverage past diet data and quantify individual food preferences as the objective function of each DM's diet optimization problem, which can generate diets that are consistent with the inferred preferences \cite{ghobadi2018robust, shahmoradi2019quantile}. 
}

\rev{
We use the proposed optimality-based clustering approach to integrate the clustering of diet decisions and the inference of objective functions representing the DMs' food preferences. We use the database from the National Health and Nutrition Examination Surveys (NHANES) that provides nutritional requirements and nutrition facts per serving to build a diet recommendation problem (see Appendix \ref{sec:Diet_Data}); to simplify the experiment, we consider 13 nutrients, classify foods into nine food ``types,''  and assume that once the number of servings for each food type is determined, decisions on specific menus can be made by dietitians post-hoc, similar to the experiment done in~\cite{shahmoradi2019quantile}. We assume that DMs approximately solve their own DMPs (i.e., the diet problems) whose constraints correspond to nutritional requirements as well as maximum serving size allowances representing food availability or budgetary restriction. While the constraint coefficients (i.e., nutrition facts) are fixed across all DMPs, the right-hand-side values---lower and upper limits on each nutrient and the maximum serving size per food type---are assumed to vary to account for different age, gender, level of physical activity, and health conditions of the DMs \cite{DietaryGuidline}. We generate 20 hypothetical DMs each with a unique DMP whose lower and upper nutrient limits and the maximum serving sizes are randomly drawn from the ranges specified in 
Appendix \ref{sec:Diet_Data}. For each DM $k$, 
we solve the corresponding DMP with an arbitrary objective function to produce a vector of food servings; to make the experiment realistic, we then add a noise from the uniform distribution [0,1] to each component of the vector. We 
treat the resulting food intake vector as an observed diet decision data point from DM $k$, denoted by $\hat\bx^{k}$. We apply the SC, CI, and IC approaches to this data to compare their performance in reproducing diets similar to the observed ones. We set $L=4$, i.e., four clusters (hence four cluster-specific objective functions) are sought. 
}

\begin{table}[h]
  \centering
  \caption{\rev{Worst-case distances for all DMs achieved by SC, CI, and IC.}}
    \begin{tabular}{>{\color{black}}l>{\color{black}}c>{\color{black}}c>{\color{black}}c>{\color{black}}c>{\color{black}}c>{\color{black}}c>{\color{black}}c>{\color{black}}c>{\color{black}}c>{\color{black}}c}
    \toprule
          & \multicolumn{10}{c}{\rev{Decision Makers}} \\
\cmidrule{2-11}    Model & 1     & 2     & 3     & 4     & 5     & 6     & 7     & 8     & 9     & 10 \\
    \midrule
    SC    & 0.859 & 0.904 & 0.962 & 0.906 & 0.948 & 0.751 & 0.835 & 0.999 & 0.986 & 0.969 \\
    CI    & 0.859 & 0.904 & 3.182 & 2.597 & 1.023 & 0.751 & 0.818 & 0.999 & 0.986 & 1.609 \\
    IC    & 0.859 & 8.382 & 0.962 & 2.391 & 8.805 & 8.594 & 5.907 & 4.588 & 10.317 & 0.969 \\
    \midrule
          & \multicolumn{10}{c}{\rev{Decision Makers}} \\
\cmidrule{2-11}    Model & 11    & 12    & 13    & 14    & 15    & 16    & 17    & 18    & 19    & 20 \\
    \midrule
    SC    & 0.923 & 0.948 & 0.970 & 0.996 & 0.992 & 0.948 & 0.943 & 0.889 & 0.523 & 0.866 \\
    CI    & 0.923 & 4.293 & 0.970 & 2.043 & 0.992 & 2.562 & 3.649 & 0.889 & 0.523 & 4.198 \\
    IC    & 4.343 & 6.145 & 0.970 & 8.708 & 2.224 & 0.948 & 8.057 & 8.167 & 8.900 & 0.776 \\
    \bottomrule
    \end{tabular}%
  \label{tab:all_results}%
\end{table}%

\rev{
Table \ref{tab:all_results} shows the worst-case distance between the observed diet and a newly generated diet for each DM $k$ (i.e., $\|\hat\bx^{k} - {\bx}^{k*}\|_{\infty}$), latter of which was generated by the objective function for the cluster that the DM was assigned to by either SC, CI, or IC 
%
(we followed the same solution procedure as described in Section \ref{sec:results_synthetic}). 
The SC results were optimal as the upper and lower bounds matched with the same clusters and objective functions.
While the distance between the observed and new diets achieved by SC did not exceed 1 across all DMs, CI and IC  both often led to diets that are far from the observed ones, with the worst-case distances of 4.293 and 10.317, respectively. The performance of CI (i.e., clustering diets first and inferring objective functions post-hoc) was better than that of IC for most patients and was comparable to that of SC for 11 out of 20 DMs.
}

\begin{figure}[h!]
	\centering
	\begin{subfigure}[t]{0.47\textwidth}
		\centering
		\includegraphics[width=3.2in]{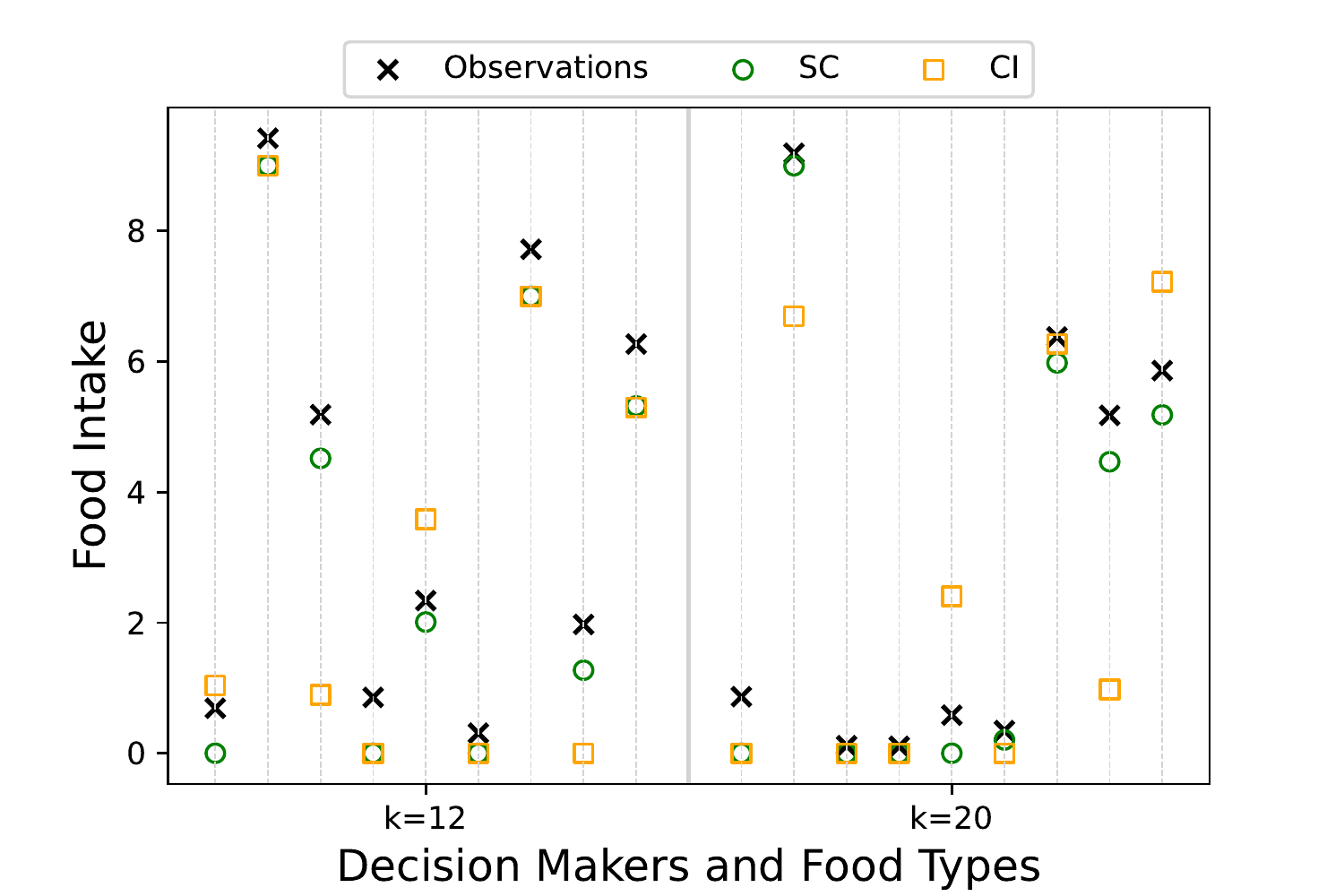}
		\caption{
		\rev{SC produces better worst-case distance than CI.}
		}
		\label{fig:Diet_High_diff}
	\end{subfigure}%
	\begin{subfigure}[t]{0.47\textwidth}
		\centering
		\includegraphics[width=3.2in]{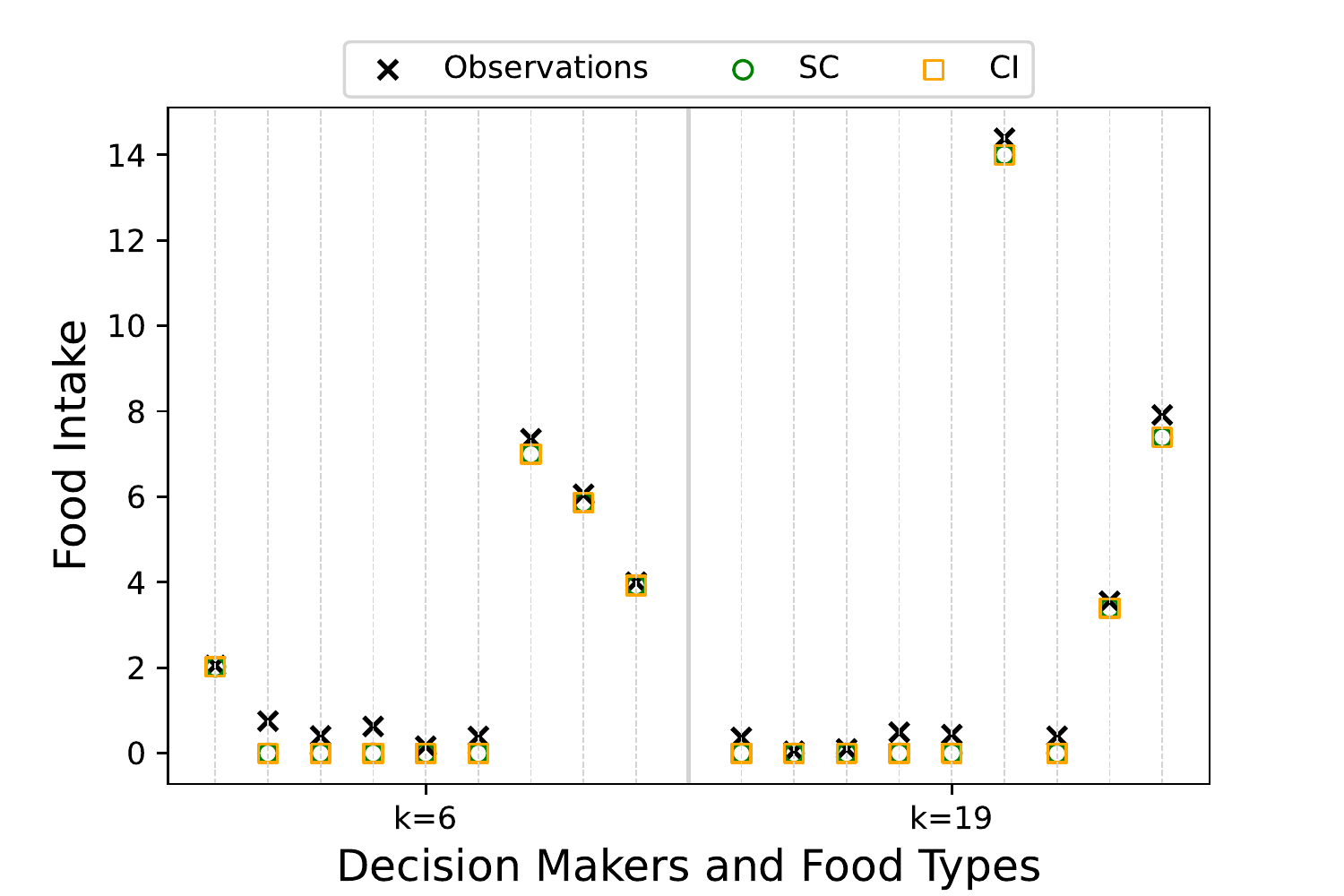}
		\caption{
		\rev{SC and CI produce same diets}
		}
		\label{fig:Diet_Low_diff}
	\end{subfigure}%
	\caption{\rev{Comparison of diet decisions generated by SC and CI for four representative DMs.}
	}
	\label{fig:Diet_Performance}
\end{figure}

\rev{
Figure~\ref{fig:Diet_Performance} shows a detailed comparison of the diet decisions generated by 
the SC and CI approaches for four representative DMs. The x-axis represents the DMs and y-axis shows the optimal servings for the nine food types for each DM. We only compare SC and CI in this figure as CI was generally better than  IC as shown in Table \ref{tab:all_results} and the CI approach can be seen similar to the existing clustering analysis in the literature in that clusters are obtained based on the similarity of past diets themselves. Figure \ref{fig:Diet_High_diff} shows two DMs ($k=12,20$) for which the SC and CI approaches led to different diets. While the diets generated by SC were generally close to the observed diets, CI generated some recommendations that were far from the observations, indicating that these DMs were assigned to clusters mixed with other DMs with different preferences and thus the resulting cluster-specific objective function generated inconsistent diets for these DMs. Figure \ref{fig:Diet_Low_diff}, on the other hand, shows that SC and CI often lead to exactly same diets. Both Table \ref{tab:all_results} and Figure~\ref{fig:Diet_Performance} indicate that the SC approach performs at least as well as CI and IC in reproducing diets consistent with the observations under varying constraints and thus the inferred cluster-specific objective function can be used for generating future diets for the individuals in the same cluster in a stable manner. A post-hoc analysis can be done to further identify association between the preference clusters and individual health conditions and socio-demographic factors.
}

\vspace{-0.1in}
\section{Conclusion}
\vspace{-0.1in}
In this paper we introduced a new clustering approach, called optimality-based clustering, that clusters DMs based on similarity of their decision preferences. We formulated the clustering problem as a non-convex optimization problem and proposed MIP formulations that provide lower and upper bounds. We also proposed two heuristics that can be efficient in large instances and perform comparably to solving the problem exactly in certain instances. \rev{We used the proposed clustering models in the diet recommendation context to cluster DMs based of their food preferences}. The future research includes extending the idea of optimality-based clustering to other types of  DMPs such as non-linear, mixed-integer, and multi-objective optimization problems.

{\small
\bibliographystyle{abbrvnat}
\bibliography{IO_Clustering}
}

\

\appendix

\noindent{\Large\textbf{Appendix}}
\vspace{-0.1in}
\section{Solution Approaches}\label{sec:solution_approaches}
While \textbf{SC-LB}($\cK,\cL$) and \textbf{SC-UB}($\cK,\cL$) 
provide bounds for the SC problem, i.e., \textbf{SC}($\cK,\cL$)
, these problems are typically large-scale MIPs and thus can be computationally challenging. To address this, we propose to use the CI and IC approaches to create initial solutions for these MIPs. However, as discussed in Section \ref{sec:sol_struct_reform}, Stage 2 of the CI approach 
itself involves $\mathbf{SC}(\cG^\ell_{\textup{CI}},L=1)$ for each $\ell\in \cL$, and Stage 1 and the post-processing step of the IC approach 
also involves solving $\mathbf{SC}(\{k\},L=1)$ for each $k\in \cK$ and $\mathbf{SC}(\cG^\ell_{\textup{IC}},L=1)$ for each $\ell\in \cL$, respectively. These ``smaller'' SC problems can also be approximately solved using the corresponding smaller upper and lower bound MIP formulations, just like how the full-size SC problem is approximated.  For example,  $\mathbf{SC}(\cG^\ell_{\textup{CI}},L=1)$ for the CI approach can be approximated by $\textbf{SC-UB}(\cG^\ell_{\textup{CI}},L=1)$ and $\textbf{SC-LB}(\cG^\ell_{\textup{CI}},L=1)$. 

\begin{figure}[h]
	\centering
	\includegraphics[width=6.4in]
	{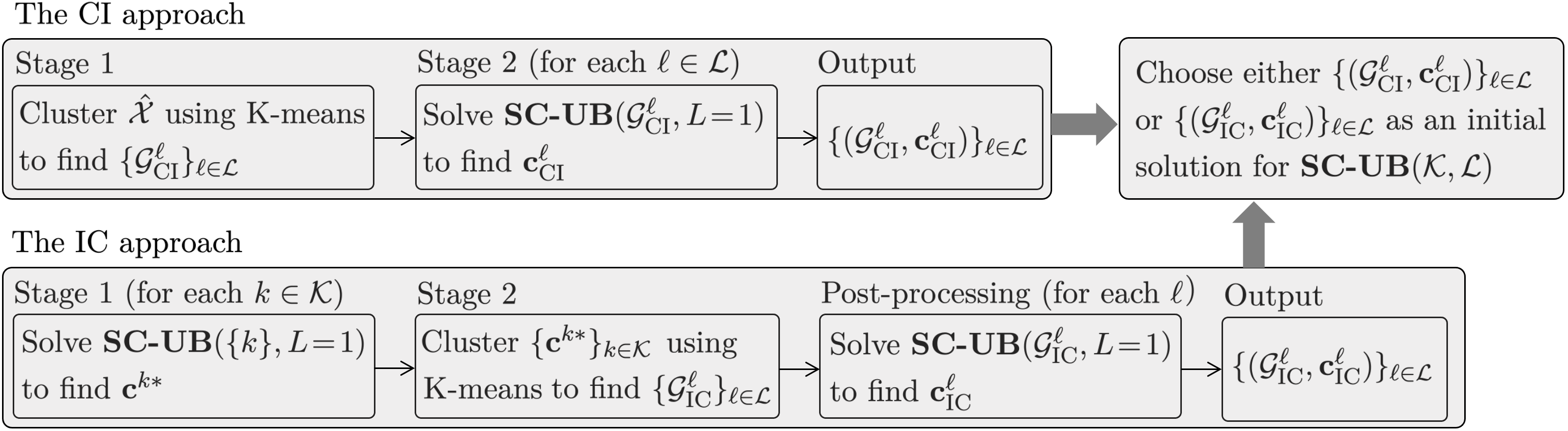}
	\caption{Steps to generate an initial solution for \textbf{SC-UB}($\cK,\cL$). Similar steps can be used for generating an initial solution for \textbf{SC-LB}($\cK,\cL$).}
	\label{fig:Solving_steps}
\end{figure}

Once approximate solutions (clusters and cost vectors) for the IC and CI approaches are obtained, we use them as  initial feasible solutions for 
\textbf{SC-LB}($\cK,\cL$) and \textbf{SC-UB}($\cK,\cL$).
In particular, solving the smaller SC problems in the CI and IC approaches using the corresponding SC-UB formulations leads to initial feasible solutions for $\textbf{SC-UB}(\cK,\cL)$; similarly, solving the smaller SC problems using the SC-LB formulations leads to initial feasible solutions for $\textbf{SC-LB}(\cK,\cL)$. Figure \ref{fig:Solving_steps} details how initial solutions are generated for $\textbf{SC-UB}(\cK,\cL)$; similar steps can be performed to find initial feasible solutions for $\textbf{SC-LB}(\cK,\cL)$. Finally, 
once $\textbf{SC-UB}(\cK,\cL)$ and $\textbf{SC-LB}(\cK,\cL)$ are solved, if their objective values are equal, their solutions will be optimal for the true SC problem. Otherwise, it is safer to use the result obtained by the upper bound formulation $\textbf{SC-UB}(\cK,\cL)$ 
because the true worst-case distance will never exceed the objective value of $\textbf{SC-UB}(\cK,\cL)$. 


\section{\rev{Data for the Diet Problem}}\label{sec:Diet_Data}

\rev{The data used for generating the diet optimization problems in Section~\ref{sec:results_diet} can be found in Table~\ref{tab:DietData}.}

\begin{threeparttable}
\caption{\zahed{Food items, nutrient data per serving, and lower and upper limits on nutrition consumption.}}
\label{tab:DietData}%
{\fontsize{8}{11} \selectfont
\begin{tabular}{lrrrrrrrrrrr}
     			\toprule
     			& \multicolumn{9}{c}{Food Type}                               &   {*Lower}    &  {*Upper} \\
     			\cmidrule{2-10}    \multicolumn{1}{c}{ } & 1  & 2 & 3 & 4  & 5 & 6  & 7  & 8 & 9 &  {Limit
			}  & {Limit
			} \\
     			\midrule
     			Energy (KCAL) & 91.53 & 68.94 & 23.51 & 65.49 & 110.88 & 83.28 & 80.50 & 63.20 & 52.16 & 1800.00 & 2500.00 \\
     			Total\_Fat (g) & 4.95  & 0.71  & 1.80  & 3.48  & 6.84  & 4.41  & 5.80  & 0.94  & 0.18  & 44.00 & 78.00 \\
     			Carbohydrate (g) & 6.89  & 12.16 & 0.25  & 0.00  & 5.44  & 4.68  & 0.56  & 11.42 & 13.59 & 220.00 & 330.00 \\
     			Protein (g) & 4.90  & 3.68  & 1.59  & 7.99  & 6.80  & 5.93  & 6.27  & 2.40  & 0.41  & 56.00 & NA \\
     			Fiber (g) & 0.00  & 0.06  & 0.00  & 0.00  & 0.28  & 0.29  & 0.00  & 1.19  & 1.81  & 20.00 & 30.00 \\
     			Vitamin C (mg) & 0.01  & 1.76  & 0.00  & 0.00  & 0.17  & 0.16  & 0.00  & 0.02  & 11.19 & 90.00 & 2000.00 \\
     			Vitamin B6 (mg) & 0.06  & 0.03  & 0.01  & 0.09  & 0.11  & 0.06  & 0.06  & 0.03  & 0.10  & 1.30  & 100.00 \\
     			Vitamin B12 (mcg) & 0.67  & 0.39  & 0.09  & 0.65  & 0.11  & 0.63  & 0.56  & 0.00  & 0.00  & 2.40  & NA \\
     			Calcium (mg) & 172.09 & 125.72 & 46.24 & 2.21  & 5.90  & 15.03 & 29.00 & 27.21 & 6.14  & 1000.00 & 2500.00 \\
     			Iron (mg) & 0.05  & 0.08  & 0.04  & 0.75  & 0.35  & 0.35  & 0.73  & 0.79  & 0.13  & 8.00  & 45.00 \\
     			Copper (mg) & 0.02  & 0.04  & 0.01  & 0.03  & 0.03  & 0.03  & 0.04  & 0.05  & 0.06  & 0.90  & 10.00 \\
     			Sodium (mg) & 61.02 & 48.24 & 65.08 & 72.32 & 211.05 & 128.27 & 223.50 & 125.62 & 1.42  & 1500.00 & 2300.00 \\
     			Vitamin A (mcg) & 42.89 & 22.24 & 12.78 & 0.00  & 1.33  & 9.53  & 81.00 & 0.01  & 13.04 & 900.00 & 3000.00 \\
     			\multirow{2}[1]{*}{**Max serving ($\pm \bar\eta$) } & \multirow{2}[1]{*}{8} & \multirow{2}[1]{*}{8} & \multirow{2}[1]{*}{8} & \multirow{2}[1]{*}{8} & \multirow{2}[1]{*}{8} & \multirow{2}[1]{*}{8} & \multirow{2}[1]{*}{8} & \multirow{2}[1]{*}{8} & \multirow{2}[1]{*}{8} & \multicolumn{2}{c}{\multirow{2}[1]{*}{}} \\
     			&       &       &       &       &       &       &       &       &       & \multicolumn{2}{c}{} \\
     			\bottomrule
\end{tabular}}
\begin{tablenotes}
	\footnotesize
	\item * Lower and upper limit values for each DM are chosen from $[(1-\eta)LL,LL]$ and $[UL, (1+\eta) UL]$, respectively, where $LL$ and $UL$ correspond to the limit values presented in the last two columns and $\eta\in [0,20]$ for each nutrient. 
	\item ** Max serving sizes are randomly chosen integers $8 \pm \bar \eta$ where $\bar\eta\in\{1,2,3,4\}$. 
  \end{tablenotes}
\end{threeparttable}

\section{Proofs}\label{sec:Proofs}

\begin{proof}{Corollary}{\ref{corollary:same_As}} 
Consider an optimal solution $\big(\{(\bc^{\ell*},\cG^{\ell*})\}_{\ell\in \cL},\{(\bx^{k*},\by^{k*})\}_{k\in \cK}\big)$ to \eqref{eq:StabClust}. Let $\cX$ denote the set of solutions for $\bx$ that satisfy $\bA\bx\ge \bbb$ (i.e., feasible for the DMP). Let $\cX^*(\bc)$ be the set of optimal solutions for the DMP with cost vector $\bc$. Since $\bA^k=\bA$ and $\bbb^k=\bbb$ for all $k\in \cK$, we have  $\cX^*(\bc)=\cX^{k*}(\bc)$ for all $k\in \cK$ and thus $\text{ext}(\cX^*(\bc))=\text{ext}(\cX^{k*}(\bc))$, $\forall k\in \cK$, for any cost vector $\bc$. Let $Q_\ell=\left|\text{ext}(\cX^*(\bc^{\ell*}))\right|$, $\cQ^\ell=\{1,\ldots,Q_\ell\}$ (i.e., $\text{ext}(\cX^*(\bc^{\ell*}))=\{\bar\bx^1,\ldots, \bar\bx^{Q_\ell}\}$), and $q^*_\ell\in \displaystyle \argmax_{q_\ell\in \cQ^\ell} \left\{ \max_{ k\in \cG^{\ell*}} \{\|\hat\bx^k-\bar\bx^{q_\ell}\|_r\}\right\}$. Let $\cI^{\ell*}=\{i\in \cI \,|\, {\ba^i}'{\bar\bx^{q^*_\ell}}=b_i\}$ for each $\ell\in \cL$. From the proof of Theorem \ref{theorem:corner_point}, we know that $\big(\{(\bc^{\ell*},\cG^{\ell*})\}_{\ell\in \cL},\{(
\bar\bx^{q^*_\ell},\by^{k*})\}_{k\in \cK}\big)$ 
is also optimal for \eqref{eq:StabClust} for each $\ell\in \cL$ and $k\in \cG^{\ell*}$. Note that because $\bar\bx^{q^*_\ell}$ is an extreme point of $\cX$, the interior of $\text{cone}(\{\ba^{i}\}_{i\in \cI^{\ell*}})$, 
{i.e., $\text{cone}_+(\{\ba^{i}\}_{i\in \cI^{\ell*}})$}, is nonempty (also see Proposition 15 in \cite{tavasli2018}). Hence, there exist $\tilde\by^\ell \ge \bzero$ and $\tilde\bc^\ell$ such that $\tilde y^\ell_i>0$ for all $i\in \cI^{\ell*}$, $\|\tilde\bc^\ell\|_1=1$, and $\tilde\bc^{\ell}=\sum_{i\in \cI^{\ell*}}\tilde y_i^{\ell} \ba^{i}$, i.e., $\tilde\bc^{\ell}\in \text{cone}_+(\{\ba^{i}\}_{i\in \cI^{\ell*}})$ for all $\ell\in \cL$. To complete the proof, given such $\tilde\by^\ell$ and $\tilde\bc^\ell$ 
we note that the solution $\big(\{(\tilde\bc^{\ell},\cG^{\ell*})\}_{\ell\in \cL},\{(
\bar\bx^{q^*_\ell},\tilde\by^{k})\}_{k\in \cK}\big)$ where $\tilde\by^k=\tilde\by^\ell$ for all $k\in \cG^{\ell*}$ is also feasible (hence optimal) for \eqref{eq:StabClust}, because $\|\tilde\bc^\ell\|_1=1$ for all $\ell\in \cL$
, $\bA\bar\bx^{q^*_\ell}\ge \bbb$, $\tilde\by^
\ell\ge\bzero$, ${\bA}'\tilde\by^
\ell=\tilde\bc^\ell$ for each $\ell\in \cL$ and $k\in \cG^{\ell*}$, and finally, ${\tilde\bc^{\ell'}}\tilde\bx^k=\sum_{i\in \cI^{\ell*}} \tilde y^k_i\ba^i\tilde\bx^k=\sum_{i\in \cI^{\ell*}} \tilde y^k_i b_i=\bbb'\tilde\by^k$ for each $\ell\in \cL$ and $k\in \cG^{\ell*}$.
\end{proof}
\begin{proof}{Proposition}{\ref{prop:strict_MIP}}
To prove part (i), we show that given an optimal solution for \eqref{eq:StabClust} we can construct a feasible solution for \eqref{eq:Re_StabClust_MIP} that achieves the objective value no greater than $\rho^*$. 
Let $\big(\{(\bc^{\ell*},\cG^{\ell*})\}_{\ell\in \cL},\{(\bx^{k*},\by^{k*})\}_{k\in \cK}\big)$ be an optimal solution to \eqref{eq:StabClust}. From the proof of Theorem \ref{theorem:corner_point}, for a fixed $\ell\in \cL$ and $k\in \cG^{\ell*}$, the set of $\bx^k$'s that satisfy \eqref{eq:StabClust_5}--\eqref{eq:StabClust_6}, together with $(\bc^{\ell*},\cG^{\ell*},\by^{k*})$, can be characterized by $\text{conv}(\text{ext}(\cX^k(\bc^{\ell*})))$ where $\text{conv}(\cdot)$ denotes the convex hull of a given set of points. Let $\cQ^k=\{1,\ldots, Q_k\}$ be the index set for the extreme points in $\text{ext}(\cX^k(\bc^{\ell*}))$. 

We now 
construct a feasible solution for \eqref{eq:Re_StabClust_MIP}. For all $\ell\in \cL$, let $\tilde c^{\ell+}_j= c^{\ell*}_j$, $\tilde c^{\ell-}_j=0$, and $\tilde z^\ell_j=1$ if $c^{\ell*}_j\ge 0$, and let $\tilde c^{\ell+}_j=0$, $\tilde c^{\ell-}_j=- c^{\ell*}_j$, and $\tilde z^\ell_j=0$ otherwise. Since $\|\bc^{\ell*}\|_1=1$, $(\bc^{\ell},\bc^{\ell+},\bc^{\ell-},\bz^\ell)=(\bc^{\ell*},\tilde\bc^{\ell+},\tilde\bc^{\ell-},\tilde\bz^\ell)$ satisfies \eqref{eq:Re_StabClust_MIP_2}--\eqref{eq:Re_StabClust_MIP_5} for each $\ell\in\cL$. Construct $\tilde \bu$ by letting $\tilde u_{k\ell}=1$ if $k\in \cG^{\ell*}$ and $\tilde u_{k\ell}=0$ otherwise, and let $\tilde \blambda^k=\by^{k*}$. 
For $k\not\in\cG^{\ell*}$, i.e., $(k,\ell)$ such that $\tilde u_{k\ell}=0$, constraint \eqref{eq:Re_StabClust_MIP_6} holds trivially with $(\blambda^k,\bu)=(\tilde\blambda^k,\tilde\bu)$ as $M_1$ is a sufficiently large positive constant. For $k\in \cG^{\ell*}$, because ${\bA^k}'\by^{k*}=\bc^{\ell*}$ from constraint~\eqref{eq:StabClust_3} we have ${\bA^k}'\tilde\blambda^{k}=\bc^{\ell*}$, which satisfies \eqref{eq:Re_StabClust_MIP_6} with $\tilde u_{k\ell}=1$. 
Also, $\tilde \lambda^k$ and $\tilde v^k_i$ satisfy \eqref{eq:Re_StabClust_MIP_7} because $\tilde \lambda^k=y^{k*}_i=0$ whenever $\tilde v_i^k=1$. 
Next, for each $k\in \cK$ and some arbitrary extreme point $\bar{q}_k \in \cQ^k$, let $\tilde\cI_{\bar{q}_k}=\{i\in \cI^k\ |\  {\ba^{ki}}'\bx^{\bar{q}_k}=b^k_i\}$; then let $\tilde v^k_i=1$ if $i\in \tilde\cI_{\bar q_k}$ and $\tilde v^k_i=0$ otherwise. Clearly, by definition of $\tilde \cI_{\bar q_k}$, we have $\bx^{\bar q_k}$ and $\tilde v_i^k$ satisfy \eqref{eq:Re_StabClust_MIP_8}. 
Furthermore, since $\bx^{\bar q_k}$ is an extreme point, we have $|\tilde \cI_{\bar q_k}|=\sum _{i\in \cI^k} \tilde v_i^k=n$, satisfying \eqref{eq:Re_StabClust_MIP_9}. Since each data point $k$ is assigned to one of the clusters $\{\cG^{\ell*}\}_{\ell\in \cL}$, we have $ \sum_{\ell\in \cL} \tilde u_{k\ell}=1$, which satisfies constraint \eqref{eq:Re_StabClust_MIP_10}. Thus, the solution $(\{(\tilde\bc^{\ell+},\tilde\bc^{\ell-},\tilde\bc^{\ell},\tilde\bz^{\ell})\}_{\ell\in \cL},\{(\tilde\blambda^k,\bx^{\bar q_k},\tilde \bv^k)\}_{k\in \cK}, \tilde\bu
 )$ is feasible for problem \eqref{eq:Re_StabClust_MIP}. Let $\tilde\beta=\displaystyle\max_{k\in \cK}\{d(\hat\bx^k,\bx^{\bar q_k})\}$, i.e., the objective function value achieved by this solution. 
 Then we have 
$\tilde\beta=\displaystyle\max_{k\in \cK}\{d(\hat\bx^k,\bx^{\bar q_k})\}\le 
\displaystyle \max_{k\in \cK} \max_{q_k\in\cQ^k}\{d(\hat\bx^k,\bx^{q_k})\}\le 
\displaystyle \max_{k\in \cK} \;
\max \{d(\hat\bx^k,\bx)\, | \,  \bx\in \text{conv}(\text{ext}(\cX^k(\bc^{\ell*})))\}=
\rho^*$. Finally, because $\beta^*\le \tilde\beta$, we have $\beta^*\le \rho^*$, as desired.
%

To prove part (ii), let $( \{\tilde\bx^k\}_{k\in \cK},\{\tilde\bv^{k}\}_{k\in \cK},\{\blambda^k\}_{k\in \cK}, 
\{\tilde \bc^{\ell}\}_{\ell\in \cL},\tilde\bu)$ be 
optimal for \eqref{eq:Re_StabClust_MIP} and assume $\tilde \bc^{\ell}\in \textup{cone}_+(\{\ba^{ki}\}_{i,k:\tilde v^{k}_i=1,\tilde u_{k\ell}=1})$. Consider $\bc^{\ell*}=\tilde \bc^\ell$, $\bx^{k*}=\tilde\bx^{k}$,  $\by^{k*}={\tilde\blambda}^k$, and $\cG^{\ell*}=\{k\in \cK \,|\, \tilde u_{k\ell}=1\}$ for all $k \in \cK$ and $\ell \in \cL$. We show that the solution $(\{\bc^{\ell*}, \cG^{\ell*}\}_{\ell\in\cL}, \{\bx^{k*}, \by^{k*}\}_{k\in\cK})$ is  feasible for \eqref{eq:StabClust} as follows. First, we have ${\bA^k}'\by^{k*}={\bA^k}'\tilde\blambda^{k}=\tilde\bc^{\ell}=\bc^{\ell*}$ for all $k \in \cG^{\ell*}$ if $\tilde u_{k\ell}=1$, which satisfies \eqref{eq:StabClust_3}; $\by^{k*}=\tilde\blambda^{k}\ge \bzero$, which satisfies \eqref{eq:StabClust_4}; $\bA^k\bx^{k*}=\bA^k\tilde\bx^k\ge \bbb^k$ for all $k\in \cK$ and $\ell\in \cL$, which satisfies \eqref{eq:StabClust_5}; and $\|\bc^{\ell*}\|_1=\|\tilde\bc^{\ell}\|_1=1$ for all $\ell\in \cL$, which satisfies constraint \eqref{eq:StabClust_2}. 
To show this solution also satisfies \eqref{eq:StabClust_6}, we let $\tilde \cI_{k}=\{i\in \cI^k\, |\, \tilde v_i^k=1\}$ for each $k\in \cK$. From \eqref{eq:Re_StabClust_MIP_8} and $\bx^{k*}=\tilde \bx^k$, we have ${\ba^{ki}} \bx^{k*}= b_i^k$ for all $i\in \tilde \cI^k$. From \eqref{eq:Re_StabClust_MIP_7} and $\by^{k*}=\tilde\blambda^k$, we have $y^{k*}_i\ge 0$ for $i\in \tilde \cI^k$ and $y^{k*}_i= 0$ otherwise. Thus, we have $\displaystyle\sum_{i\in \cI^k}y^{k*}_i{\ba^{ki}} \bx^{k*}= \sum_{i\in \cI^k}y^{k*}_i b_i^k$, and because $\displaystyle\sum_{i\in \cI^k}y^{k*}_i{\ba^{ki}}=\bc^{\ell*}$ from \eqref{eq:Re_StabClust_MIP_6}, this equation becomes ${\bc^{\ell*}}'\bx^{k*}={\bbb^k}'\by^{k*}$, which satisfies \eqref{eq:StabClust_6}. As a result, $(\{\bc^{\ell*}, \cG^{\ell*}\}_{\ell\in\cL}, \{\bx^{k*}, \by^{k*}\}_{k\in\cK})$ is  feasible for \eqref{eq:StabClust}. 
Next, note that $\bc^{\ell*}\in \text{cone}_+(\{\ba^{ki}\}_{i\in \tilde \cI^k})$ and $|\tilde \cI^k|=n$ for each $k\in \cG^{\ell*}$; therefore, $\bx^{k*}$ is an extreme point and in fact is the only solution for $\bx^k$ that satisfies constraints in \eqref{eq:StabClust} for each $k\in\cK$, i.e., $\text{ext}(\cX^{k*}(\bc^{\ell*}))=\{\bx^{k*}\}$ and $|\cQ^k|=1$ for all $\ell\in \cL$ and $k\in \cG^{\ell*}$. Thus, we have 
$\displaystyle \tilde\beta 
=  \max_{\ell\in \cL} \max_{k: \tilde u_{k\ell}=1} d(\hat\bx^k,\tilde\bx^{k}) 
= \max_{\ell\in \cL}\max_{k \in \cG^{\ell*}} d(\hat\bx^k,\tilde\bx^{k}) = \max_{\ell\in \cL} \max_{k \in \cG^{\ell*}} d(\hat\bx^k,\bx^{k*}) = \max_{\ell\in \cL} \max_{k \in \cG^{\ell*}}\max_{q_k\in \cQ^k}\{d(\hat\bx^k,\bx^{q_k})\}= \rho^*,$ where the fourth equality holds because $\bx^{k*}$ is the only member in $\cQ^k$.
\end{proof}
\begin{proof}{Proposition}{\ref{prop:strict_MIP_UB}}
Let $(\{\bx^{k*}, \bv^{k*}, \blambda^{k*}\}_{k\in \cK}, \{ \bc^{\ell*}\}_{\ell\in \cL}, \bu^{*})$ be optimal for problem \eqref{eq:Re_StabClust_MIP_UB} and  $\beta^*$ be the optimal value of \eqref{eq:Re_StabClust_MIP_UB}. Due to constraint \eqref{eq:Re_StabClust_MIP_UB_3}, we have $\bc^{\ell*}\in \textup{cone}_+(\{\ba^{ki}\}_{i,k:{v}^{k*}_i=1, u_{k\ell}^*=1})$ for each $\ell\in \cL$. Note that this solution is feasible for \eqref{eq:Re_StabClust_MIP} with the same objective value $\beta^*$ because the feasible region of \eqref{eq:Re_StabClust_MIP_UB} is a subset of that of \eqref{eq:Re_StabClust_MIP} due to the extra constraint \eqref{eq:Re_StabClust_MIP_UB_3}. 
From the proof of Proposition~\ref{prop:strict_MIP}~(ii), if this solution is optimal for \eqref{eq:Re_StabClust_MIP}, then $\beta^*$ is equal to the optimal value of \eqref{eq:StabClust}, i.e., $\beta^*=\rho^*$; on the other hand, if this solution is feasible for \eqref{eq:Re_StabClust_MIP}, $\beta^*\ge\rho^*$. Thus, $\beta^*$ is an upper bound on $\rho^*$.
\end{proof}

\end{document}